\documentclass[10pt,landscape,amsart]{article}
\usepackage{amsmath,amsbsy,amsfonts}
\title{\sc Stabilit\'e en niveau 0, pour les
 groupes orthogonaux
impairs p-adiques.}
\date{}
\author{Moeglin, Colette\\
\sl CNRS, Institut de math\'ematiques de JUSSIEU\\
moeglin@math.jussieu.fr\\ }
\begin{document}
\maketitle
\vskip 0.5cm Pr\'ecisons tout de suite que dans ce qui suit,
$F$ est un corps extension finie de $\mathbb{Q}_p$ avec
$p\neq 2$ et m\^eme pour le th\'eor\`eme principal $p$ grand.
Le but de ce travail est de produire des fonctions sur les
groupes p-adiques orthogonaux impairs dont les int\'egrales
orbitales sur les \'el\'ements elliptiques r\'eguliers ne
d\'ependent que des  classes de conjugaison stable. Au
passage, on  produit aussi des fonctions dont la somme des
int\'egrales orbitales
\`a l'int\'erieur d'une classe stable fix\'ee est nulle. A la
fin du papier, on interpr\`ete ce r\'esultat en terme de
stabilit\'e des repr\'esentations elliptiques de niveau 0
pour ces groupes orthogonaux. Si l'on a bien pr\'edit les
signes qui d\'ependent encore de \cite{auber}, c'est la somme
des repr\'esentations dans un paquet qui est stable et ces
combinaisons lin\'eaires engendrent l'espace distributions
stables combinaison lin\'eaires de repr\'esentations
elliptiques.

Dans le d\'etail, on commence par rappeler ce qu'est une
repr\'esentation de niveau z\'ero et comment on peut  lui
associer un pseudo-coefficient; ce n'est pas  nouveau et
n'est pas au c{\oe}ur du papier. C'est quand on passe \`a la
description des param\`etres de ces repr\'esentations que
l'on entre dans le vif du sujet. On d\'ecrit ces param\`etres
en terme d'ensemble d'orbites unipotentes de groupes complexes
convenables et de syst\`emes locaux sur ces orbites; la fa\c
con classique de faire cela est de consid\'erer le morphisme
de Langlands-Lusztig, $\psi$, de $W_F\times SL(2,{\mathbb C})$
dans
$Sp(2n,{\mathbb C})$ (le groupe dual) et de d\'ecomposer
d'abord la restriction de $\psi$ \`a $W_F$. Ici on
d\'ecompose la restriction de $\psi$ au sous-groupe de $W_F$
noyau de l'application $W_F \rightarrow <Fr>$ (o\`u
$Fr$ est un Frob\'enius) et c'est dans le commutant de
l'image par $\psi$ de ce sous-groupe que vivent les orbites
unipotentes et syst\`emes locaux ci-dessus. Pour faire cette
classification, on utilise le fait que
$\psi$ est de niveau
$0$ (cf la d\'efinition donn\'ee dans le texte) mais ceci
n'est pas une hypoth\`ese importante \`a cet endroit. On peut
alors utiliser les m\'ethodes de Lusztig et la
repr\'esentation de Springer g\'en\'eralis\'ee pour associer
aux syst\`emes locaux trouv\'es, des fonctions sur des groupes
finis, qui sont en fait des parahoriques en r\'eduction du
groupe orthogonal de d\'epart; c'est l\`a que l'hypoth\`ese
de niveau 0 intervient. Ce sont les faisceaux caract\`eres de
Lusztig. On remonte ces fonctions sur le groupe parahorique
en les rendant invariantes par le radical pro-p-unipotent et
on les prolonge par z\'ero pour en faire des fonctions sur le
groupe orthogonal. Ce sont ces fonctions pour lesquelles on
peut calculer le comportement des int\'egrales sur les
classes de conjugaison d'\'el\'ements elliptiques \`a
l'int\'erieur d'une classe de conjugaison stable. Avec cette
m\'ethode, ce ne sont pas des combinaisons lin\'eraires de
param\`etres de Langlands qui donnent des objets stables mais
directement certains param\`etres; c'\'etait d\'ej\`a le cas
en \cite{mw} et ici c'est expliqu\'e en \ref{stabilite}. On
revient \`a des combinaisons lin\'eraires plus habituelles en
faisant une op\'eration style transformation de Fourier, cf.
\ref{fourier}.

Pour finir, on veut interpr\'eter 
ces fonctions comme un ensemble de pseudo coefficients pour 
les repr\'esentations elliptiques (elliptiques au sens
d'Arthur) de niveau z\'ero (cf. \ref{conjecture}); vu ce qui
est rappel\'e au d\'ebut de ce papier, pour le faire on doit
calculer ce que l'on appelle la restriction aux parahoriques
des repr\'esentations, c'est-\`a-dire calculer l'action de
chaque parahorique dans la sous-repr\'esentation form\'ee par
les vecteurs invariants sous l'action du radical
pro-p-unipotent du dit parahorique. Pour cela, on a besoin de
2 r\'esultats. Le premier est un r\'esultat annonc\'e par
Aubert, Kutzko et Morris (\cite{auber}) qui ram\`ene l'\'etude
de l'alg\`ebre de Hecke des repr\'esentations induites de
cuspidales de niveau 0 pour un groupe r\'eductif \`a des
alg\`ebres de Hecke de repr\'esentations induites \`a partir
de cuspidales de r\'eduction unipotente pour des groupes
r\'eductifs convenables; ici ce r\'esultat fera intervenir
d'autres groupes orthogonaux et des groupes unitaires et il y
aura sans doute un signe. Il faut
ensuite savoir calculer la restriction aux parahoriques des
repr\'esentations de r\'eduction unipotente pour ces groupes
orthogonaux et  unitaires. Pour les groupes orthogonaux, c'est
essentiellement fait en
\cite{w2} (nous l'avons d\'ej\`a utilis\'e dans \cite{mw}) et
pour les groupes unitaires c'est fait dans \cite{m}.
Moyennant le travail en cours \cite{auber} on a donc une
bonne description de pseudo coefficient pour les
repr\'esentations elliptiques de niveau 0 des groupes
orthogonaux consid\'er\'es ici. Et en utilisant les
r\'esultats d'Arthur (\cite{arthur}) qui ram\`enent les
probl\`emes de stabilit\'e pour des repr\'esentations
elliptiques \`a la stabilit\'e des int\'egrales orbitales en
les \'el\'ements elliptiques de leurs pseudo-coefficients, on
en d\'eduit une description
 des paquets stables de repr\'esentations
elliptiques de niveau z\'ero. Tout ceci est r\'eminiscent de
\cite{mw}.

Pour finir cette introduction, je remercie Anne-Marie Aubert
pour les conversations que nous avons eues et le texte
qu'elle a \'ecrit pour moi, ainsi que Jean-Loup Waldspurger
qui m'a fait un nombre certain de calculs.

\section{Repr\'esentations de niveau z\'ero.}

En suivant les d\'efinitions usuelles et en particulier
celles de  \cite{mw},  on appelle r\'eseau presque autodual un
r\'eseau $L$ de $V$ tel que
\[\omega_F\tilde{L} \subset L\subset \tilde{L}.
\] Pour $L$ un r\'eseau presque autodual comme ci-dessus, on
note $K(L)$ le stabilisateur du r\'eseau et $U(L)$ le radical
pro-p-unipotent.  Plus g\'en\'eralement, on appelle
cha\^{\i}ne de r\'eseaux presque autoduaux une famille
$L_.=(L_0,L_1, \cdots, L_r)$ de r\'eseaux de $V$ telle que:
\[\omega_F\tilde{L}_r \subset L_r\subset \cdots
\subset L_1 \subset L_0\subset
\tilde{L}_0.\] Et on g\'en\'eralise de fa\c con \'evidente la
d\'efinition de $K(L_.)$ et $U(L_.)$. Une description
totalement explicite de  ces objets a \'et\'e donn\'ee en 
\cite{mw} 1.2. On y a en particulier d\'efini la notion
d'association.

Soit $\pi$ une repr\'esentation irr\'eductible de $G(F)$; on
dit que $\pi$ est de niveau 0 s'il existe un r\'eseau presque
autodual $L$ tel que l'espace des invariants sous
$U(L)$, $\pi^{U(L)}$ est non nul; c'est une d\'efinition
standard reprise en particulier de  \cite{auber}. On note
alors
$\pi_L$ la repr\'esentation de $K(L)/U(L)$ dans ces
invariants; c'est une repr\'esentation (non 
irr\'eductible en g\'en\'eral) d'un groupe r\'eductif sur le
corps
${\mathbb F}_q$.

La notion de repr\'esentation de niveau z\'ero pour les
groupes $GL$ est de m\^eme ordre et nous l'utiliserons.

La construction de pseudo-coefficients pour les s\'eries
discr\`etes (ou plus g\'en\'eralement les repr\'esentations
elliptiques au sens d'Arthur) faite en  \cite{mw} 1.9 (iv)
s'\'etend au cadre des repr\'esentations elliptiques de
niveau 0. C'est ce que nous allons expliquer dans cette
partie.

\subsection{Support cuspidal\label{support}}  Soit $\pi$ une
repr\'esentation irr\'eductible de $G(F)$; on appelle support
cuspidal de $\pi$ la donn\'ee d'un sous-groupe de Levi $M$ de
$G$, qui est le Levi d'un parabolique d\'efini sur $F$ et une
repr\'esentation cuspidale $\pi_{cusp}$ de $M(F)$ tel que
$\pi$ soit quotient de l'induite de $\pi_{cusp}$ (gr\^ace \`a
un parabolique de Levi $M$). La donn\'ee de $(M,\pi_{cusp})$
est alors d\'efinie \`a conjugaison pr\`es par le groupe de
Weyl de $G$.

A partir de maintenant, on suppose que $\pi$ est de niveau
z\'ero. Il r\'esulte de \cite{mp} 6.11 que 
$\pi_{cusp}$ l'est aussi (le groupe $G$ \'etant remplac\'e
par $M$). On v\'erifie alors que l'on peut construire  une
cha\^{\i}ne de r\'eseaux presque autoduaux, $L_.$ de $V$, en
bonne position par rapport \`a $M$, telle que
$\pi_{cusp}$ ait des invariants (non nuls) sous $U(L_.)\cap
M(F)=:U_M(L_.)$; on note $\pi_{cusp}^{U_M(L_.)}$ cet espace
d'invariants. On note
$K_M(L_.):=K(L_.)\cap M(F)$ et
$\pi_{cusp}^{U_M(L_.)}$ est naturellement une
repr\'esentation du groupe fini $K_M(L_.)/U_M(L_.)$; elle est
cuspidale. On note $\chi_{cusp}$ la donn\'ee de cette
repr\'esentation et du groupe fini qui op\`ere; cela
sous-entend que la donn\'ee d'une cha\^{\i}ne de r\'eseaux
presque autoduaux, $L_{.,cusp}$ ait \'et\'e faite. Un tel 
choix n'est pas unique mais il l'est \`a association pr\`es.

Il r\'esulte de \cite{mp} que pour $L$ un r\'eseau presque
autodual de $V$, la repr\'esentation $\pi_L$ n'est pas nulle
si et seulement si $K(L)$ est associ\'e \`a un sous-groupe
parahorique contenant $K(L_{.,cusp})$ et le support cuspidal
de $\pi_L$ comme repr\'esentation du groupe fini $K(L)/U(L)$
est conjugu\'ee de $\chi_{cusp}$. On note
$C(K(L))_{\chi_{cusp}}$ l'ensemble des fonctions sur le
groupe fini $K(L)/U(L)$ engendr\'e par les caract\`eres des
repr\'esentations irr\'eductible ayant un conjugu\'e de
$\chi_{cusp}$ comme support cuspidal.   Et on note
$C(K(L))_{cusp,\chi_{cusp}}$ la projection de l'espace
$C(K(L))_{\chi_{cusp}} $ sur l'ensemble des fonctions
cuspidales.

Lusztig associe aux repr\'esentations des
groupes finis et donc \`a $\chi_{cusp}$ un
\'el\'ement semi-simple $s_{\chi}$ dans un certain groupe sur
$\overline{\mathbb{F}}_q$, le groupe
$Sp(2n',\overline{\mathbb{F}}_q)\times
O(2n'',\overline{\mathbb{F}}_q)$. La classe de
conjugaison de
$s_{\chi}$ dans
$GL(2n,\overline{\mathbb{F}}_q)$ est elle bien d\'efinie. On
dira que $s_\chi$, ou plut\^ot sa classe de conjugaison, est
la classe de conjugaison semi-simple associ\'ee \`a
$\chi_{cusp}$.

\subsection{Repr\'esentation elliptique\label{elliptique}}
Comme en  \cite{mw} 1.7, on utilise la notion de
repr\'esentation elliptique telle que pr\'ecis\'ee par
Arthur; c'est, une repr\'esentation elliptique est une
combinaison lin\'eaire de repr\'esentations temp\'er\'ees. Si
l'on fixe un Levi $M$ de $G$, Levi d'un sous-groupe
parabolique de $G$ d\'efini sur $F$ et une repr\'esentation
cuspidale $\pi_{cusp}$ de $M(F)$, on dit que la
repr\'esentation elliptique, $\pi$, est de support cuspidal
$(M,\pi_{cusp})$ si toutes les repr\'esentations
irr\'eductibles qui interviennent dans sa combinaison
lin\'eaire ont cette propri\'et\'e; et on dit qu'elle est de
niveau 0 si le support cuspidal est de niveau 0.

\subsection{Pseudo-coefficients\label{pseudo}} On reprend les
notations de \ref{support}, en particulier
$\chi_{cusp}$ et $C(K(L))_{cusp,\chi_{cusp}}$. En suivant
\cite{mw}, on remonte tout \'el\'ement $f$ de
$C(K(L))_{cusp,\chi_{cusp}}$, en une fonction sur $K(L)$
invariante par $U(L)$ et on la prolonge en une fonction
not\'ee $f^G$ sur $G(F)$ en l'\'etendant par 0 hors de
$K(L)$. Par cette proc\'edure, on obtient une fonction
cuspidale sur $G(F)$, c'est-\`a-dire une fonction dont les
int\'egrales orbitales sur les \'el\'ements semi-simples non
elliptiques sont nulles. Quand on somme cette construction
sur tous les supports cuspidaux de niveau 0, on construit
ainsi un morphisme de $\oplus_{L/\sim}C(K(L))_{cusp}$ dans
l'ensemble des fonctions cuspidale sur $G(F)$.

On note $Ps^G_{ell,\chi_{cusp}}:=\oplus_{L/\sim}
C(K(L))_{cusp,\chi_{cusp}}$ que l'on munit d'un produit
scalaire. Pour  d\'efinir ce produit scalaire, il faut fixer
un ensemble de repr\'esentants des groupes $K(L)$ modulo
association (les sommes sur $K$ ci-dessous signifient la
somme sur un ensemble de repr\'esentants)
\[ (\sum_{K}\phi_{K},\sum_{K}\phi'_{K})=\sum_K
w(K)^{-1}(\phi_K,\phi'_K)_K,\] o\`u le dernier produit
scalaire est le produit scalaire usuel sur un groupe fini et
o\`u $w(K)$ est un volume d\'ecrit en \cite{mw} 1.6.
Quand on somme cette construction sur l'ensemble des supports
cuspidaux de niveau 0, on d\'efinit $Ps^G_{ell,0}$ muni d'un
produit scalaire. La d\'ecomposition suivant les supports
cuspidaux (modulo conjugaison) est une somme directe.

On rappelle la construction des pseudo coefficients pour les
repr\'esentations elliptiques de niveau 0 faite
(essentiellement) en \cite{mw} 1.9. et qui repose sur les
travaux de \cite{ss} et de \cite{courtes}

Pour tout $K$ comme ci-dessus, d\'efinissons $B(K)$ (resp.
$B(K)_{\chi_{cusp}})$ une base de $C(K)_{cusp}$ (resp.
$C(K)_{cusp,\chi_{cusp}}$) et pour toute repr\'esentation
virtuelle temp\'er\'ee, $D$, posons 
\[\phi_D:=\oplus_K \sum_{f\in B(K)}w(K) \overline{D(f^G)} f;
\] 
\[\phi_{D,\chi_{cusp}}:=\oplus_K \sum_{f\in
B(K)_{\chi_{cusp}}}w(K)
\overline{D(f^G)} f.
\]  {\bf Th\'eor\`eme: }{\sl L'application $D\mapsto \phi_D$
induit un isomorphisme de l'espace engendr\'e par les
caract\`eres des repr\'esentations elliptiques de
niveau 0 sur
$Ps^G_{ell,0}$. Cet isomorphisme est compatible \`a la
d\'ecomposition suivant le support cuspidal.} 

\rm On peut r\'ecrire $\phi_D$ sous la forme la plus
utilisable. On pose:
\[ \varphi_D:=\sum_K w(K) tr_K(D),\] o\`u $tr_K(D)$ est la
trace pour la repr\'esentation de
$K/U_K$ dans l'espace des invariants de la repr\'esentation
$D$ sous $U_K$ (le radical pro p-unipotent de $K$). On note
$proj_{ell}\varphi_D$ la projection de $\varphi_D$ sur les
fonctions cuspidales (cette projection se fait pour chaque
parahorique $K$ individuellement). On sait alors d\'efinir
$(proj_{ell}\varphi_D)^G$ qui est un pseudo coefficient de
$D$ si $D$ est elliptique (cf
\cite{mw} 1.9 (iv)). Si
$D$ a $\chi_{cusp}$ (cf. ci-dessus) pour support cuspidal,
alors $proj_{ell}\varphi_D\in \oplus_K C[K]_{cusp}$

{\bf Remarque: }{\sl Avec les notations ci-dessus,
$(proj_{ell}\varphi_D)^G$ est un pseudo-coefficient de la
repr\'esentation elliptique $D$.}

On rappelle aussi que d'apr\`es Arthur \cite{arthur} 6.1, 6.2
(on a enlev\'e l'hypoth\`ese relative au lemme fondamental en
\cite{mw} 4.6) une combinaison lin\'eaire $D$ de
repr\'esentations elliptiques est stable si et seulement si
les int\'egrales orbitales de $\varphi_D^G$ sont constantes
sur les classes de conjugaison stable d'\'el\'ements
elliptiques r\'eguliers.

\section {Classification des param\`etres discrets de
niveau 0
}

 On consid\`ere les couples
$(\psi,\epsilon)$ de morphismes continus suivants:
\[
\begin{array}{c}\psi: W_F\times SL(2,{\mathbb C})  \rightarrow
Sp(2n,{\mathbb C})\\
\epsilon: Cent_{Sp(2n,{\mathbb C})}(\psi) \rightarrow
\{\pm 1\},\end{array}
\] o\`u, en notant $I_F$ le sous-groupe de ramification
 de $W_F$, la restriction de $\psi$ \`a $I_F$ est triviale sur
le groupe de ramification sauvage et o\`u le centralisateur
de $\psi$ dans $Sp(2n,{\mathbb C})$ n'est pas inclus dans un
sous-groupe de Levi de $Sp(2n,{\mathbb C})$. De tels couples
sont appel\'es des param\`etres discrets de niveau 0; comme
la condition ne porte que sur $\psi$ et non sur $\epsilon$,
on peut dire aussi que $\psi$ est discret de niveau 0 sans
r\'ef\'erence \`a $\epsilon$.

\subsection {Morphismes de
param\'etrisation\label{classification}}

 Pour donner la classification des morphismes comme ci-dessus,
il est plus simple d'avoir fix\'e un g\'en\'erateur du groupe
ab\'elien
$I_F/P_F$, o\`u $P_F$ est le groupe de ramification sauvage.
Et pour cela, il est plus simple de fixer une extension
galoisienne mod\'er\'ement ramifi\'ee $E/F$ et de ne
consid\'erer que les morphismes $\psi$ qui se factorisent par
le groupe de Weil relatif $W_{E/F}$. On fixe $Fr$ une
image r\'eciproque d'un Frob\'enius de l'extension non
ramifi\'ee dans $W_{E/F}$ et $s_E$ un g\'en\'erateur du
groupe multiplicatif du corps r\'esiduel de $E$. Dans ce cas
la restriction de $\psi$ \`a $I_F$ est d\'etermin\'ee par
l'image de $s_E$. A conjugaison pr\`es c'est donc la donn\'ee
des valeurs propres de la matrice image de $s_E$ par $\psi$
qui d\'etermine cette restriction. On va donc fixer cette
restriction en la notant $\chi$, c'est \`a dire fixer une
matrice de $Sp(2n,{\mathbb C})$ dont les valeurs propres sont
des racines de l'unit\'e d'ordre premier \`a
$p$. On peut donc oublier $E$ et garder $\chi$ et
consid\'erer que $\chi$ est d\'etermin\'e par une collection
de racines de l'unit\'e, l'ensemble des valeurs propres
ensemble que l'on note $VP(\chi)$. Pour $u\in VP(\chi)$ on
note $mult(u)$ la multiplicit\'e de $u$ en tant que valeur
propre. L'action du Frobenius transforme $\chi$ en $\chi^q$,
ainsi si $u\in VP(\chi)$ alors $u^q\in VP(\chi)$ et
$mult(u)=mult(u^q)$. Comme $\chi$ est \`a valeurs dans
$Sp(2n,{\mathbb C})$, l'espace propre pour la valeur propre
$u$ est en dualit\'e avec l'espace propre pour la valeur
propre $u^{-1}$, d'o\`u aussi $mult(u)=mult(u^{-1})$.  A
l'int\'erieur de $VP(\chi)$ on d\'efinit l'\'equivalence
engendr\'ee par la relation \'el\'ementaire
$u\sim u^q$. On note $[VP(\chi)]$ les classes d'\'equivalence
et si $u\in VP(\chi)$, on note $[u]$ sa classe
d'\'equivalence. On v\'erifie que s'il existe
$u\in VP(\chi)$ tel que $u^{-1}\notin [u]$ alors le
centralisateur de $\psi$ est inclus dans un sous-groupe de
Levi de $Sp(2n,{\mathbb C})$; on suppose donc que $u^{-1}\in
[u]$ pour tout $u$. Pour tout
$u\in VP(\chi), u\notin \{\pm 1\}$, on d\'efinit $\ell_{[u]}$
comme le plus petit entier tel que $u^{-1}=u^{q^{\ell_[u]}}$;
le cardinal de la classe
$[u]$ est alors $2\ell_{[u]}$. On pose:
$m([u]):=mult(u)$ o\`u $u \in VP(\chi)$ dans la classe de
$[u]$ comme la notation le sugg\`ere. On remarque que 
\[2n=m(1)+m(-1)+\sum_{[u]\in [VP(\chi)],u\notin\{\pm
1\}}m([u])2\ell_{[u]} .\] Soit $M$ un entier et $U$ une
orbite unipotente de
$GL(M,{\mathbb C})$; on dit que $U$ est symplectique (resp.
orthogonale)
 si tous les blocs de Jordan sont pairs (resp. impairs) et on
dit qu'elle est discr\`ete si son nombre de blocs de Jordan
d'une taille donn\'ee est au plus 1. D'o\`u la notation
symplectique discr\`ete et orthogonale discr\`ete qui allie
les 2 d\'efinitions.

\bf Proposition: \sl l'ensemble des homomorphismes $\psi$
ci-dessus (c'est-\`a-dire discrets et de niveau 0), pris 
\`a conjugaison pr\`es, dont la restriction \`a $I_F$ est
conjugu\'ee de $\chi$  est en bijection avec l'ensemble des
collections d'orbites unipotentes $\{U_{[u],\zeta}, [u] \in
[VP(\chi)], \zeta\in \{\pm 1\}\}$, de groupe
$GL(m([u],\zeta),{\mathbb C})$, ce qui d\'efinit l'entier
$m([u],\zeta)$ (\'eventuellement 0) avec les propri\'et\'es
suivantes:
\[\forall{[u]\in [VP(\chi)]}, \, m([u],+)+m([u],-)=m([u]);\]
pour tout $[u]\in [VP(\chi)], u\neq \pm 1$, l'orbite
$U_{[u],+}$ est une orbite symplectique discr\`ete, l'orbite
$U_{[u],-}$ est une orbite orthogonale discr\`ete et  les
orbites $U_{[\pm 1],\pm }$ sont des orbites symplectiques
discr\`etes.\rm

L'int\'er\^et de ramener la classification \`a une collection
d'orbites unipotentes est de pouvoir ensuite utiliser la
repr\'esentation de Springer g\'en\'eralis\'ee pour
construire des repr\'esentations de groupes de Weyl, puisque
l'on aura aussi des syst\`emes locaux sur ces orbites.

On a d\'ecrit avant l'\'enonc\'e comment on comprenait la
restriction de $\psi$ \`a $I_F$; pour avoir la restriction de
$\psi$ \`a $W_F$, il faut encore d\'ecrire l'image du
rel\`evement du Frob\'enius, $Fr$, \`a conjugaison pr\`es.
Par commodit\'e et uniquement dans cette d\'emonstration, on
note $V$ l'espace vectoriel ${\mathbb C}^{2n}$ et pour
$u\in 	VP(\chi)$, on note $V[u]$ l'espace propre
correspondant \`a cette valeur propre. Les conditions que
doivent v\'erifier $\psi(Fr)$ sont: \^etre une matrice
symplectique et induire un isomorphise entre $V[u]$ et
$V[u^q]$ pour tout $u\in VP(\chi)$.

Pour traduire ces conditions, fixons $u\in VP(\chi)$. Il faut
distinguer les 2 cas:

premier cas: $u\neq \pm 1$. On remarque que
$\psi(Fr^{q^{2\ell_u}})$ induit un isomorphisme de $V[u]$
dans lui-m\^eme. On note
$F_u$ cet homomorphisme. Le groupe
$GL(V[u])$ s'identifie naturellement \`a un sous-groupe de
$Sp(2n,{\mathbb C})$. Comme on ne cherche \`a classifier les
morphismes 
$\psi$ qu'\`a conjugaison pr\`es, on peut encore conjuguer
sous l'action de $GL(V[u])$; cela se traduit sur
$F_u$ par la conjugaison habituelle. A conjugaison pr\`es
$F_u$ est donc d\'etermin\'e par ses valeurs propres dont on
note
$VP(F_u)$ l'ensemble. On v\'erifie encore que si $VP(F_u)$
contient un
\'el\'ement autre que $\pm 1$, alors l'image de $\psi$ est
incluse dans un sous-groupe de Levi propre de $Sp(2n,{\mathbb
C})$. Pour
$\zeta\in \{\pm 1\}$, on note $V[u,\zeta]$ l'espace propre
pour la valeur propre $\zeta$ de $F_u$. On remarque pour la
suite que $V[u,\zeta]$ est muni du produit scalaire:
$$\forall v,v' \in V[u,\zeta], <v,v'>_u:=<v,
\psi(Fr)^{q^{\ell_u}}v'>.$$ Et, pour $v$ et $v'$ comme
ci-dessus:
$$<v,\psi(Fr)^{q^{\ell_u}}v'>=<\psi(Fr)^{-q^{\ell_u}}v,v'>=
\zeta<\psi(Fr)^{-q^{\ell_u}} F_u
v,v'>=\zeta<\psi(Fr)^{q^{\ell_u}}v,v'>$$
$$=-\zeta <v',\psi(Fr)^{q^{\ell_u}}v>=-\zeta<v',v>_u.
$$ Ainsi, la forme $<\, ,\, >_u$ est symplectique pour
$\zeta=1$ et orthogonale pour $\zeta=-1$. Il est facile de
v\'erifier que cette forme est non d\'eg\'en\'er\'ee. Ces
constructions se font donc pour tout $u\in VP(\chi)$
diff\'erent de $\pm 1$. De plus 
$\psi(Fr)$ induit une isom\'etrie de
$V[u,\zeta]$ sur $V[u^q,\zeta]$; ceci permet de d\'efinir
intrins\`equement l'espace orthogonal ou symplectique
$V([u],\zeta)$ pour tout $[u]\in [VP(\chi)]$ muni du produit
scalaire $<\, ,\, >_{[u]}$.

deuxi\`eme cas: $u\in \{\pm 1\}$. On d\'efinit ici $F_u$
comme l'action de $\psi(Fr)$ comme automorphisme de $V[u]$.
On v\'erifie comme ci-dessus que si $F_u$ a des valeurs
propres autres que $\pm 1$, l'image de $\psi$ se trouve dans
un Levi de $Sp(2n,{\mathbb C})$; on d\'efinit donc encore
$V[u,\zeta]$ pour $\zeta=\pm 1$ les valeurs propres de
$F_u$. Mais ces espaces sont ici des espaces symplectiques par
restriction de la forme symplectique.

Comme
$\psi(SL(2,{\mathbb C}))$ commute \`a $\psi(W_F)$ les images
des
\'el\'ements unipotents de $SL(2,{\mathbb C})$ s'identifient
\`a des
\'el\'ements unipotents des automorphismes des espaces
$V([u],\zeta)$ pour tout $[u]\in [VP(\chi)]$ et tout $\zeta
\in \{\pm 1\}$. A conjugaison pr\`es, le morphisme $\psi$
restreint \`a $SL(2,{\mathbb C})$ est m\^eme uniquement
d\'etermin\'e par l'orbite de ces \'el\'ements. Ce sont ces
orbites qui sont not\'ees $U_{[u],\zeta}$ dans l'\'enonc\'e.
Comme les \'el\'ements de $\psi(SL(2,{\mathbb C}))$ commutent
\`a  $\psi(Fr)$ et respectent la forme symplectique, ils
respectent chaque forme $<\, ,\, >_u$. Ce sont donc des
orbites unipotentes du groupe d'automorphismes de la forme.
Il reste \`a remarquer que si l'une de ces orbites a 2 blocs
de Jordan de m\^eme taille, alors l'image de $\psi$ est
incluse dans un Levi. R\'eciproquement la donn\'ee des
orbites permet de reconstruire (\`a conjugaison pr\`es)
l'homomorphisme
$\psi$.

\bf Remarque: \sl Soit $\chi$ comme ci-dessus et identifions
les racines de l'unit\'e d'ordre premier \`a $p$ de ${\mathbb
C}$ avec leurs analogues dans
$\overline{\mathbb{F}}_q$. Les \'el\'ements de $VP(\chi)$
avec leur multiplicit\'e d\'efinissent donc un \'el\'ement de
$GL(2n,\overline{\mathbb{F}}_q)$ dont la classe de
conjugaison est bien d\'efinie.

\rm Avec cette remarque, on peut associer \`a $\chi$ un
\'el\'ement semi-simple $s_\chi$ bien d\'efini \`a
conjugaison pr\`es dans $GL(2n,\overline{\mathbb{F}}_q)$.
C'est l'analogue du $s_{\chi}$ de \ref{support}.

\subsection{Syst\`eme local  \label{centralisateur}} On fixe
$\psi,\epsilon$ commme dans l'introduction de cette section et
on reprend les notations de la preuve pr\'ec\'edente en notant
$Jord(U_{[u],\zeta})$, o\`u $[u]\in [VP(\chi)]$ et
$\zeta\{\pm 1\}$, l'ensemble des blocs de Jordan des
orbites unipotentes associ\'ees \`a $\psi$.

\bf Remarque: \sl le centralisateur de $\psi$ est isomorphe
\`a $\prod_{[u]\in VP(\chi)}\prod_{\alpha\in
Jord(U_{[u],\zeta})}\{\pm 1\}$. L'image du centre de
$Sp(2n,{\mathbb C})$ dans ce commutant est l'\'el\'ement
$-1$ diagonal. Ainsi $\epsilon$ s'identifie \`a une
application de $\cup_{[u]\in [VP(\chi)];\zeta\in \{\pm
1\}}Jord(U_{[u],\zeta})$ dans $\{\pm 1\}$.

\rm
En reprenant la preuve pr\'ec\'edente, on voit que le
commutant de $\psi$ s'identifie au commutant de
$\psi(SL(2,{\mathbb C}))$ vu comme sous-ensemble de
$\times_{[u]\in
[VP(\chi)]}\times_{\zeta=\pm}Aut(V([u],\zeta),<\, ,\, >_u)$.
On sait calculer ce commutant.  C'est alors un produit de
groupes orthogonaux
$\times_{[u],\zeta}\times_{\alpha\in
Jord(U_{[u],\zeta})}O(mult_\alpha,{\mathbb C})$, o\`u
$mult_\alpha$ est la multiplicit\'e de $\alpha$ comme bloc de
Jordan de l'orbite en question; pour $\psi$ discret cette
multiplicit\'e est 1. Pour avoir ce r\'esultat la seule
hypoth\`ese utilis\'ee est que $U_{[u],\zeta}$ est
symplectique si
$<\, ,\, >_u$ est symplectique et orthogonale sinon.  On aura
aussi
\`a regarder le cas elliptique  o\`u cette hypoth\`ese sur le
type de $U_{[u],\zeta}$ est satisfaite mais pas la
multiplicit\'e 1; on utilisera alors cette description. Dans
le cas de la multiplicit\'e 1, le groupe orthogonal se
r\'eduit \`a $\{\pm 1\}$; d'o\`u l'\'enonc\'e,
l'identification du centre \'etant imm\'ediate.

Remarquons encore que quelle que soit la multiplicit\'e, on
peut voir le
$\epsilon$ comme une application de
$\times_{[u],\zeta}Jord(U_{[u],\zeta})$ dans $\{\pm 1\}$.

\section{Faisceaux caract\`eres.}

\subsection{Construction de fonctions.\label{faisceaux}} 

Soit $m\in \mathbb{N}$; on utilisera fr\'equemment la
notation $D(m)$ pour l'ensemble des couples d'entiers
$(m',m'')$ tels que $m=m'+m''$. On fixe
$\chi$ un morphisme comme en \ref{classification}.  On
reprend les notations
$[VP(\chi)]$ de \ref{classification}. Pour tout $[u]\in
[VP(\chi)]$ avec $[u]\neq \pm 1$, on a d\'efini les entiers
$m([u])$ (qui sont les multiplicit\'es des valeurs propres).
On pose $n([u])=m([u])$ si $[u]\neq \pm 1$ et
$n(1)=m(1)/2$, $n(-1)=m(-1)/2$.  Pour
${(n'_{[u]},n''_{[u]})\in D(n([u])}$, on note ${\mathbb
C}[\hat{W}_{n'_{[u]},n''_{[u]}}]:= {\mathbb
C}[\hat{\mathfrak{S}}_{n'{[u]}}]\otimes {\mathbb
C}[\hat{\mathfrak{S}}_{n''{[u]}}] $ et
${\mathbb C}[\hat{W}_{D([u]}]$ l'espace vectoriel
$\oplus_{(n'_{[u]},n''_{[u]})\in D(n([u])}
{\mathbb
C}[\hat{W}_{n'_{[u]},n''_{[u]}}] $, o\`u les
chapeaux repr\'esentent les classes d'isomorphie de
repr\'esentations du groupe chapeaut\'e. Pour $u=\pm 1$, la
situation est plus compliqu\'ee \`a cause de l'existence de
faisceaux caract\`eres cuspidaux.  On garde la m\^eme
notation (pour unifier) mais on remplace
$\hat{\mathfrak{S}}_{n'[u]}$ et
$\hat{\mathfrak{S}}_{n''[u]}$  par l'ensemble des symboles de
rang $n'[u]$ respectivement $n''[u]$ de d\'efaut impair
respectivement pair; il est rappel\'e en
\cite{w3} 2.2, 2.3 comment ces symboles param\'etrisent aussi
des repr\'esentations irr\'eductibles de groupes; un symbole
de d\'efaut impair, $I=:2h+1$, et de rang $n'([u])$
param\'etrise une repr\'esentation du groupe de Weyl de type
$C$ et de rang $n'([u])-h^2-h$.  Dans
 le cas du d\'efaut pair, il faut admettre les d\'efauts
n\'egatifs; dans la r\'ef\'erence donn\'ee tout est expliqu\'e
avec pr\'ecision, les difficult\'es venant de la non
connexit\'e des groupes orthogonaux pairs et du fait que pour
un tel groupe il faut regarder simultan\'ement la forme
d\'eploy\'ee et celle qui ne l'est pas. Grosso modo, un
symbole de d\'efaut pair, $2h''$, et de rang $n''([u])$
param\'etrise une repr\'esentation d'un groupe de Weyl de
type C de rang $n''([u])-(h'')^2$.

Fixons maintenant un ensemble de paires
$(n'{[u]},n''{[u]})\in D(n([u])$. On pose:
$$ n':=\sum_{[u]\in [VP(\chi)]; [u]\neq [\pm
1]}n'{[u]}\ell_{[u]}  \, + (n'{[1]}+n'{[-1]}),$$
$$ n'':=\sum_{[u]\in [VP(\chi)]; [u]\neq [\pm
1]}n''{[u]}\ell_{[u]}\, + (n''{[1]}+n''{[-1]}).
$$ On pose $\sharp=iso$ si $G$ est d\'eploy\'e et $\sharp=an$
sinon. On note alors ${K}_{n',n''}$ un sous-groupe parahorique
(non connexe) de $G$ dont le groupe en r\'eduction, 
$\overline{K}_{n',n''}$ est isomorphe \`a $SO(2n'+1,
\mathbb {F}_q)\times O(2n'',\mathbb{F}_q)_{\sharp}$ (cf.
\ref{support}). Il est bien d\'efini \`a association pr\`es. 
 On note $M$ un sous-groupe de $\overline{K}_{n',n''}$
isomorphe \`a
$$\times_{[u]\neq [\pm 1]} U(n'{[u]},{\mathbb
{F}}_{q^{2\ell{[u]}}}/\mathbb{F}_{q^{\ell_u}})\times
SO(2(n'[1]+n'[-1])+1,\mathbb{F}_q)\qquad$$
$$ \qquad
\times_{[u]\neq [\pm 1]} U(n''{[u]},{\mathbb
{F}}_{q^{2\ell_{[u]}}}/\mathbb{F}_{q^{\ell_u}})\times
O(2(n''[1]+n''[-1]),\mathbb{F}_q)_{\sharp};$$ ci-dessus, on
n'a pas pr\'ecis\'e le plongement car cela n'a pas
d'importance, sur les corps finis il n'y a qu'une classe
de formes unitaires. Gr\^ace 
\`a Lusztig (\'etendu au cas non connexe cf.
\cite{w3} 3.1 et 3.2), on sait associer \`a un
\'el\'ement de ${\mathbb
C}[\hat{W}_{\underline{n'},\underline{n''}}]$ et \`a $\chi$
une fonction sur M, la trace du faisceau caract\`ere
associ\'e. Puis on d\'efinit cette fonction sur
$\overline{K}_{n',n''}$ (par induction); c'est une fonction
invariante par conjugaison. 

En sommant sur toutes les d\'ecompositions $D(\chi)$, on
construit ainsi une application de ${\mathbb
C}[\hat{W}_{D(\chi)}]$ dans l'ensemble des fonctions 
$\oplus_{n',n''\in D(n)}{\mathbb C}[\overline{K}_{n',n''}]$.
On remonte ensuite de telles  fonctions en des fonctions sur
$K_{n',n''}$ par invariance et on les prolonge \`a
$SO(2n+1,F)_{\sharp}$ par 0.  On note
$k_{\sharp,\chi}$ cette application. Quand on fait une somme
directe de
$\sharp=iso$ avec $\sharp=an$, on la note $k_{\chi}$.

\subsection{Support cuspidal des faisceaux
caract\`eres\label{bis}}

Dans cette section, on fixe quelques notations relatives aux
faisceaux caract\`eres  quadratiques unipotents; elles
viennent essentiellement (\`a des modifications formelles
pr\`es) de
\cite{w3} 3.1 et 3.7  lui-m\^eme fortement inspir\'e de
Lusztig. Les difficult\'es viennent de la pr\'esence de
faisceaux caract\`eres cuspidaux; c'est le cas des groupes
orthogonaux impairs et pairs qui a \'et\'e sommairement
exp\'edi\'e ci-dessus qu'il faut pr\'eciser. 

Pour les groupes orthogonaux impairs,
$SO(2m'+1,\mathbb{F}_q)$ on forme les faisceaux caract\`eres
quadratiques unipotents avec la donn\'ee d'un couple
ordonn\'e de 2 symboles  de d\'efaut impair dont la somme des
rangs est
$m'$. Notons
$\Lambda'_{+},\Lambda'_-$ ces 2 symboles et $I_{+},I_-$ leurs
d\'efauts. On \'ecrit encore $I_{\pm}=:2h'_{\pm}+1$ en
utilisant le fait que les d\'efauts sont impairs. On retrouve
alors \`a peu pr\`es les notations de \cite{w3} 3.7. On
consid\`ere le couple d'entier $(h'_++h'_-+1, \vert
h'_+-h'_-\vert)$ et le signe $\sigma':=+$ si $h'_+\geq h'_-$
et $-$ sinon. Dans ce couple d'entiers, l'un des nombres est
pair et l'autre est impair; on note $r'_p$ celui qui est pair
et $r'_{im}$ celui qui est impair. On note $n'_{\pm }$ le
rang de $\Lambda'_{\pm}$ et on pose
$N'_{\pm}:=n'_{\pm}-h'_{\pm}(h'_{\pm}+1)$. Ainsi
$\Lambda'_{\pm}$ param\'etrise une repr\'esentation du groupe
de Weyl de type $C$ de rang $N'_{\pm}$. Tandis que le couple
$r'_{im},\sigma' r'_p$ d\'etermine un faisceau cuspidal pour
le groupe $SO(r_{im}^{'2}+r_p^{'2},\mathbb{F}_q)$ et l'on a:
$2N'_++2N'_-+r_{im}^{'2}+r_p^{'2}=2m'+1$.

Pour les groupes orthogonaux pairs,
$O(2m'',\mathbb{F}_q)_{\sharp}$, on forme un faisceau
caract\`ere quadratique unipotent \`a l'aide d'un couple
ordonn\'e de 2 symboles eux-m\^emes ordonn\'es au sens qu'un
symbole est form\'e de 2 ensembles de nombres (avec des
propri\'et\'es). Au sens habituel, l'ordre des ensembles n'a
pas d'importance et le d\'efaut est la diff\'erence entre le
cardinal de l'ensemble ayant le plus d'\'el\'ement (au sens
large) et celui de l'ensemble ayant le moins d'\'el\'ements
(au sens large). Ici, les 2 ensembles sont ordonn\'es et le
d\'efaut est la diff\'erence entre le cardinal du premier
ensemble et celui du deuxi\`eme, ainsi le d\'efaut
peut-\^etre n\'egatif. On demande uniquement que les
d\'efauts soient pairs (0 est un nombre pair). On note
$\Lambda''_+,\Lambda''_-$ le couple des 2 symboles et
$P_+,P_-$ la valeur absolue de leur d\'efaut et
$\zeta_+,\zeta_-$ les signes des d\'efauts; on fera une
convention sur le signe quand le d\'efaut est 0 ci-dessous,
pour le moment on n'en a pas besoin. Ainsi
$P_\pm$ sont des nombres positifs ou nuls pairs.  On pose
encore $r''_\pm:=(\zeta_+P_+\pm \zeta_-P_-)/2$; on a ainsi 2
\'el\'ements de $\mathbb{Z}$ de m\^eme parit\'e. On note
$n''_{\pm}$ le rang de $\Lambda''_{\pm}$  et
$N''_{\pm}:=n''_{\pm}-(h''_{\pm})^2$ (o\`u $h''_{\pm}=1/2
P_{\pm}$).  Ainsi
$\Lambda''_{\pm}$ param\'etrise une repr\'esentation du groupe
de Weyl de type $C$ de rang $N''_{\pm}$. Tandis que le couple
$r''_{+}, r''_-$ d\'etermine un faisceau cuspidal pour le
groupe
$O(r_{+}^{''2}+r_-^{''2},\mathbb{F}_q)$ (cf. \cite{w3} 3.1) et
l'on a:
$2N'_++2N'_-+r_{+}^{''2}+r_-^{''2}=2m''$.

On aura \`a consid\'erer simultan\'ement 2 couples ordonn\'es
form\'e chacun de 2 symboles
$(\Lambda'_\epsilon,\Lambda''_\epsilon);\epsilon\in \{\pm\}$
o\`u, pour $\epsilon=+$ ou $-$, $\Lambda'_\epsilon$ est  de
d\'efaut impair,
$I_\epsilon$ et $\Lambda''_\epsilon$ est de d\'efaut pair
$\zeta_\epsilon P_\epsilon$ avec $P_\epsilon\in
\mathbb{N}$ et $\zeta_\epsilon\in \{\pm \}$ avec ici la
convention que si $P_\epsilon=0$ alors
$\zeta_\epsilon=(-1)^{(I_\epsilon-1)/2}$. 

Pour $\epsilon=+1$ ou $-1$, on pose pr\'ecis\'ement
$\hat{W}_{D(n[\epsilon ])}$ l'ensemble des couples de symboles
$\Lambda'_{\epsilon},\Lambda''_\epsilon$ comme ci-dessus dont
la somme des rangs vaut $n[\epsilon ]$. Ainsi
$\hat{W}_{D(n[+1])}\times \hat{W}_{D(n[-1])}$ est un ensemble
en bijection avec l'ensemble des quadruplets de symboles
ordonn\'es dont le  premier et le troisi\`eme  sont de
d\'efaut impair et les 2 autres de d\'efaut pair avec des
conditions sur la somme des rangs. On pourra donc
interpr\'eter cet ensemble en utilisant ce qui est ci-dessus
comme un ensemble des couples de repr\'esentations
quadratiques unipotentes des groupes
$SO(2m'+1,\mathbb{F}_q)\times O(2m'',\mathbb{F}_q)$ o\`u
$m'+m''=n[+1]+n[-1]$. Avec cette
interpr\'etation et ce que l'on a vu ci-dessus, les d\'efauts
des symboles d\'eterminent des faisceaux cuspidaux,
c'est-\`a-dire, combinatoirement, des nombres entiers
$\underline{r}:=(r'_+,r''_+,r'_-,r''_-)$ et un signe $\sigma'$
avec
$r'_+$ positif et impair,
$r'_-$ positif ou nul et pair et $r''_+,r''_-$ des entiers
relatifs de m\^eme parit\'e. On pose alors $\vert\underline{
r}\vert $ le quadruplet $( r'_+,\vert r''_+\vert,r'_-,\vert
r''_-\vert)$. C'est lui qui permet de construire des
fonctions de Green utiles pour la localisation (cf.
\ref{localisation}).

On note ${\mathbb C}[\hat{W}_{D(n[+1])}]\otimes {\mathbb
C}[\hat{W}_{D(n[-1])}]$ l'espace vectoriel de base
$\hat{W}_{D(n[+1]}\times \hat{W}_{D(n[-1]}$.

\subsection{Repr\'esentation de
Springer-Lusztig\label{springerlusztig}}

On fixe $(\psi,\epsilon)$ un param\`etre discret de niveau 0
et on note encore $\chi$ la restriction de $\psi$ au groupe
de ramification de $W_F$. A un tel param\`etre, on a
associ\'e une collection d'orbites $U_{[u],\zeta}$ o\`u
$[u]\in [VP(\chi)]$ et $\zeta\in \{\pm 1\}$ et $\epsilon$
s'identifie \`a un caract\`ere du groupe des composantes du
centralisateur d'un \'el\'ement de $U_{[u],\zeta}$; on voit
donc $\epsilon$ comme un morphisme de
$\cup_{[u],\zeta}Jord(U_{[u],\zeta})$ dans $\{\pm 1\}$ (cf.
\ref{centralisateur}). 

Pour
$[u]\in [VP(\chi)], [u]\neq \pm 1$, on pose
$U_{[u]}:=U_{[u],+}\cup U_{[u],-}$ ou plut\^ot l'orbite
unipotente de $GL(m([u]),{\mathbb C})$ engendr\'ee et on
pose:=
$$ n'{[u]}_{\psi,\epsilon}:=\sum_{\alpha\in Jord(U_{[u]});
\epsilon(\alpha)=+1}\alpha,\qquad
n''[u]_{\psi,\epsilon}:=\sum_{\alpha\in Jord(U_{[u]});
\epsilon(\alpha)=-1}\alpha.
$$ On d\'efinit alors $U'_{[u]}$ comme l'orbite unipotente de
$GL(n'([u])_{\psi,\epsilon},{\mathbb C})$ ayant comme bloc de
Jordan l'ensemble des $\alpha$ blocs de Jordan de $U_{[u]}$
pour lesquels $\epsilon(\alpha)=+$. On d\'efinit de m\^eme
$U''_{[u]}$.

  Pour
$u=\pm 1$, on pose:
$$n'[u]_{\psi,\epsilon}=\sum_{\alpha\in
Jord(U_{[u],+}}\alpha, \qquad
n''[u]_{\psi,\epsilon}=\sum_{\alpha\in Jord(U_{[u],-}}\alpha.
$$ Pour unifier les notations, on pose ici aussi
$U'_{[u]}:=U_{[u],+}$ et $U''_{[u]}:=U_{[u],-}$. Cette
collection de paires
$(n'[u]_{\psi,\epsilon},n''[u]_{\psi,\epsilon})$ est
naturellement not\'ee
$\underline{n}'_{\psi,\epsilon},\underline{n''}
_{\psi,\epsilon}$ et on voit la repr\'esentation de
Springer-Lusztig comme l'\'el\'ement de ${\mathbb
C}[\hat{W}_{\underline{n}'_{\psi,\epsilon},\underline{n''}
_{\psi,\epsilon}}]$ d\'efini ainsi:

soit $[u]\in [VP(\chi)], [u]\neq \pm 1$; Springer a associ\'e
\`a l'orbite $U'_{[u]}$ une repr\'esentation de
$\mathfrak{S}_{n'[u]_{\psi,\epsilon}}$, non irr\'eductible en
g\'en\'eral, dans la cohomologie de la vari\'et\'e des Borel
(on regarde toute la repr\'esentation pas seulement celle en
degr\'e maximal). Cela d\'efinit donc un
\'el\'ement de ${\mathbb
C}[\hat{\mathfrak{S}}_{n'[u]_{\psi,\epsilon}}]$. On fait la
m\^eme construction en rempla\c cant $U'_{[u]}$ par 
$U''_{[u]}$ et on obtient un \'el\'ement de ${\mathbb C}
[\hat{\mathfrak{S}}_{n''[u]_{\psi,\epsilon}}]$.

Soit maintenant $u=\pm 1$. Ce sont les constructions de
Lusztig qui sont rappel\'ees en \cite{mw} 5.5 (et \cite{w2}
5.1). Ici la situation est un peu plus compliqu\'ee puisque
l'on a 4 orbites les $U_{u,\epsilon'}$, pour
$u,\epsilon'\in \{\pm 1\}$ avec des syst\`emes locaux
et non pas 2 comme dans \cite{mw}. A chacune de ces orbites,
$U_{u,\epsilon'}$ avec son syst\`eme local est
associ\'e par la correspondance de Springer
g\'en\'eralis\'ee, un entier not\'e $k_{u,\epsilon'}$
et une repr\'esentation non irr\'eductible en g\'en\'eral du
groupe de Weyl de type $C$,
$W_{N_{u,\epsilon'}}$, o\`u l'on a pos\'e
$N_{u,\epsilon'}:=1/2(\sum_{\alpha\in
Jord(U_{u,\epsilon'})}\alpha -k_{u,\epsilon'}(k_{u
,\epsilon'}+1))$. On consid\`ere les 2 couples index\'es par
le choix d'un
\'el\'ement $u$ dans $\{\pm 1\}$
$(k_{u,+}+k_{u,-}+1, \vert
k_{u,+}-k_{u,-}\vert)$ et les 2 signes
$\zeta_{u}$ qui sont le signe de
$k_{u,+}-k_{u,-}$ quand  ce nombre est non
nul; s'il est nul le signe est $(-1)^{k_{u,+}}$ par
convention. Dans les couples l'un des nombres est impair et
on le note
$I_{u}$ et l'autre est pair et est not\'e
$P_{u}$. En regardant le produit tensoriel de la
repr\'esentation de $W_{N_{u,+}}$ avec celle de
$W_{N_{u,-}}$, on obtient une repr\'esentation du
produit qui se traduit en terme de couples de symboles dont
le premier est de d\'efaut $I_u$ et le deuxi\`eme de
d\'efaut $\zeta_u P_u$. C'est donc ainsi que
l'on construit un \'el\'ement de $\mathbb{C}[\hat{{\cal
W}}_{n'(u)_{\psi,\epsilon},n''(u)_{\psi,\epsilon}}]$.

\subsection{Induction,
restriction\label{inductionrestriction}}

Je ne connais pas d'autre justification aux constructions
faites ci-dessous que le fait que le r\'esultat
\'enonc\'e en \cite{mw} 5.5 et d\'emontr\'e en \cite{w2}
sugg\`ere la conjecture de \ref{conjecture}. 

Il s'agit de construire une application $\rho\circ \iota$ de
${\mathbb C}[\hat {W}_{D(\chi)}]$ dans lui-m\^eme. Cela
provient d'un produit tensoriel d'applications $\rho_{[u]}
\circ \iota_{[u]}$ de m\^eme nature pour tout $[u]\in
[VP(\chi)]$. Ces applications sont d\'efinies en \cite{mw}
3.18 pour $[u]=[\pm 1]$ et \cite{mw}  3.1 et 3.2 dans le cas
de $[u]\neq [\pm 1]$; on en rappelle  la d\'efinition
d'autant que l'on en donne une pr\'esentation un peu
diff\'erente.

Consid\'erons le cas o\`u $[u]\neq \pm 1$; on note ${
W}_{m[u]}$ le groupe de Weyl de type $C$ et de rang
$m([u])$. Pour $(n'[u],n''[u])\in D(m[u])$, on d\'efinit de
m\^eme
${ W}_{n'[u]}, { W}_{n''[u]}$; il existe une application
naturelle de
${ W}_{n'[u]} \times { W}_{n''[u]}$ sur
$\mathfrak{S}_{n'[u]}\times
\mathfrak{S}_{n''[u]}$. On peut ainsi remonter des
repr\'esentations de $\mathfrak{S}_{n'[u]}\times
\mathfrak{S}_{n''[u]}$ en des repr\'esentations de 
${ W}_{n'[u]} \times { W}_{n''[u]}$; ensuite on tensorise la
repr\'esentation obtenue par le caract\`ere $sgn_{CD}$ de
${ W}_{n''[u]}$. Puis on induit pour trouver un
\'el\'ement de ${\mathbb C}[\hat{{ W}}_{m[u]}]$.
L'application
$\iota_{[u]}$ est la somme sur toutes les paires dans
$D(m[u])$ de toutes ces op\'erations;
$\iota_{[u]}$ d\'efinit alors un isomorphisme de
$$\oplus_{(n'[u],n''[u])\in D(m[u])} {\mathbb
C}[\hat{\mathfrak{S}}_{n'[u]}\times
\hat{\mathfrak{S}}_{n''[u]}] \rightarrow {\mathbb C}[\hat{{
W}}_{m[u]}].$$

Fixons encore $(n'[u],n''[u])\in D(m[u])$. On voit maintenant
$\mathfrak{S}_{n'[u]}\times
\mathfrak{S}_{n''[u]}$ comme un sous-ensemble de ${\cal
W}_{n'[u]} \times {\cal W}_{n''[u]}$. Il y a en fait 2 fa\c
cons presque naturelles d'envoyer le groupe $\mathfrak{S}_m$
dans le groupe ${\cal W}_m$ ($m\in \mathbb N$); la premi\`ere
est l'homomorphisme \'evident $\sigma \mapsto w$ avec $w(\pm
i)=\pm \sigma(i)$ pour tout $i\in [1,m]$. La deuxi\`eme fa\c
con n'est pas un homomorphisme de groupe car elle est
d\'efinie par $\sigma\mapsto w$ avec
$w(\pm i)=\mp
\sigma(i)$; bien que cette application n'est pas un morphisme
de groupe, elle est
\'equivariante pour l'action adjointe. En revenant
\`a notre inclusion cherch\'ee c'est le produit de la
premi\`ere fa\c con appliqu\'ee
\`a $\mathfrak{S}_{n'[u]}$ avec la deuxi\`eme appliqu\'ee \`a
$\mathfrak{S}_{n''[u]}$. Cela permet alors de restreindre des
\'el\'ements de ${\mathbb C}[\hat {{ W}}_{m[u]}]$ en des
\'el\'ements de ${\mathbb C}[\hat{\mathfrak{S}}_{n'[u]}\times
\hat{
\mathfrak{S}}_{n''[u]}]$. En sommant ces constructions sur
toutes les paires dans $D(m[u])$, on obtient $\rho_{[u]}$.
Contrairement \`a $\iota_{[u]}$, $\rho_{[u]}$ n'est pas un
isomorphisme mais ce qui est important mais qui n'intervient
que de fa\c con cach\'ee dans \ref{conjecture} est que le
compos\'e
$\rho_{[u]}\circ
\iota_{[u]}$ est un isomorphisme si on se limite aux
fonctions \`a support dans les \'el\'ements
$U$-elliptiques, c'est-\`a-dire aux permutations qui se
d\'ecomposent en produit de cycles de longueur impaire (cf.
loc. cit.).

On va d\'ecrire d'une autre fa\c con cette application
$\rho\circ \iota$ pr\'ecis\'ement quand on se limite aux
permutations qui se d\'ecomposent en produit de cycles de
longueur impaire, en utilisant le fait qu'induire puis
restreindre peut aussi se faire en sens inverse, d'abord
restreindre puis induire. Pour cela soit
$m\in
\mathbb{N}$; on note
$D(m)$ l'ensemble des couples $(m',m'')$ tels que
$m=m'+m''$ et $DD(m)$ l'ensembles des quadruplets
$(m^{i,j};i,j\in \{',''\})$ tels que $\sum_{i,j}m^{i,j}=m$.
Soit
$(m',m'')
\in D(m)$ et
$(m^{i,j})\in DD(m)$; on dit que
$(m^{i,j})<_d (m',m'')$ si $m'=m^{','}+m^{',''}$ et
$(m^{i,j})<_e (m',m'')$ si $m'=m^{','}+m^{'','}$ ($d$ est
pour direct et $e$ pour entrelac\'e). On notera plus
g\'en\'eriquement
$\underline{m}$ et $\underline{\underline{m}}$ les
\'el\'ements de $D(m)$ et de $DD(m)$. Pour
$\underline{m}\in D(m)$, on consid\`ere de fa\c con
\'evidente le groupe
$\mathfrak{S}_{\underline{m}}=\mathfrak{S}_{m'}\times
\mathfrak{S}_{m''}$; on d\'efinit de m\^eme  les groupes
$\mathfrak{S}_{\underline{\underline{m}}}$ pour
$\underline{\underline{m}}\in DD(m)$. On note
$\chi_{\underline{\underline{m}}}$ le signe
$(-1)^{m^{'',''}}$.

Fixons $\underline{m} \in D(m)$; pour
$\underline{\underline{m}}\in D(m)$, on d\'efinit
l'application
$res_{d,\underline{m},\underline{\underline{m}}}$ de
${\mathbb C}[\hat{\mathfrak{S}}_{\underline{m}}]$ dans
${\mathbb C}[\hat{\mathfrak{S}}_{\underline{\underline{m}}}]$
comme l'application de restriction \'evidente si
$\underline{\underline{m}}<_d \underline{m}$ et $0$ sinon. Et
on note $ind_{e,\underline{m},\underline{\underline{m}}}$
l'application de ${\mathbb
C}[\hat{\mathfrak{S}}_{\underline{\underline{m}}}]$ dans
${\mathbb C}[\hat{\mathfrak{S}}_{\underline{m}}]$ qui est
l'induction si $\underline{\underline{m}}<_e \underline{m}$
et $0$ sinon; ici il faut, pour l'induction, consid\'erer
l'inclusion naturelle de
$\mathfrak{S}_{\underline{\underline{m}}}$ dans
$\mathfrak{S}_{\underline{m}}$ o\`u l'on \'echange d'abord 2e
et 3e facteur.

\bf Remarque: \sl Fixons $\underline{m}_0\in D(m)$ et 
consid\'erons
 $\rho\circ \iota$ comme une application de
${\mathbb C}[\hat{\mathfrak{S}}_{\underline{m}_0}]$ dans
$\oplus_{\underline{m}\in D(m)}{\mathbb
C}[\hat{\mathfrak{S}}_{\underline{m}}]$. On a alors, en se
limitant aux fonctions invariantes de support l'ensemble des
permutations ayant des cycles de longueur impaire:
$$\rho\circ \iota =\oplus_{\underline{m}\in D(m)}
\sum_{\underline{\underline{m}}\in DD(m)}
ind_{e,\underline{m},\underline{\underline{m}}} \circ
\biggl(\chi_{\underline{\underline{m}}} 
res_{d,\underline{m}_0,\underline{\underline{m}}}\biggr)$$\rm
Le seul point est de remarquer que sur les permutations
n'ayant que des cycles de longueur impaire le signe
$\chi_{\underline{\underline{m}}}$ co\"{\i}ncide avec
$sgn_{CD}$ tel que d\'efinit ci-dessus. 

\bf D\'efinition: \sl dans la suite, on d\'efinit $\rho\circ
\iota$ comme dans la remarque ci-dessus.
\rm

\medskip
Les d\'efinitions du cas $[u]=[\pm 1]$ sont plus
compliqu\'ees (cf. \cite{mw} en particulier 3.18 et 3.19) \`a
cause de la partie cuspidale. On les pr\'esente ainsi. On fixe
$\epsilon\in
\{\pm 1\}$; on doit d\'efinir une application de
$\hat{W}_{D(n[\epsilon])}$ dans
$\mathbb{C}[
\hat{W}_{D(n[\epsilon]})]$ que l'on prolongera lin\'eairement
en un endomorphisme de $\mathbb{C}[
\hat{W}_{D(n[\epsilon]})]$. Et on veut l'interpr\'eter comme
une restriction suivie d'une induction tordue  par un
caract\`ere. Un \'el\'ement de
$\hat{W}_{D(n[\epsilon ])}$ est la donn\'ee de deux symboles
l'un de d\'efaut impair et l'autre de d\'efaut pair dont la
somme des rangs est
$n'[\epsilon ]$. De fa\c con beaucoup plus compliqu\'ee  mais
\'equivalente (et qui permet de parler de repr\'esentations)
c'est la donn\'ee:

d'un entier impair $I_\epsilon$, d'un entier pair
$P_\epsilon$, d'un signe $\zeta_\epsilon$, de 2 entiers
$N'_\epsilon,N''_\epsilon$ tous ces nombres v\'erifiant
l'\'egalit\'e $N'_\epsilon+N''_\epsilon+
(I_\epsilon^2+P^2_\epsilon-1)/4=n[\epsilon]$ et de 2
repr\'esentations irr\'eductibles l'une de
$W_{N'_{\epsilon}}$ et l'autre de $W_{N''_{\epsilon}}$, les
groupes de Weyl de type $C$ et de rang \'ecrit en indice. 
On  suppose que
$\zeta_\epsilon=(-1)^{(I_\epsilon-1)/2}$ si
$P_\epsilon=0$.

On pose $\tilde{\zeta}_\epsilon:=(-1)^{(I_\epsilon 
-1)/2}\zeta_\epsilon$; en particulier, on a
$\tilde{\zeta} =1$ si $P_\epsilon=0$.

On note $\tilde{\chi}$ le caract\`ere trivial si
$I_\epsilon>P_\epsilon$ et le caract\`ere $sgn_{CD}$ sinon.

On d\'efinit une application de
$\mathbb{C}[\hat{W}_{{N}'_\epsilon}\times
\hat{W}_{{N}''_\epsilon}]$ dans
$\oplus_{(M'_\epsilon,M''_\epsilon)\vert
M'_\epsilon+M''_\epsilon=N'_{\epsilon}+N''_\epsilon}
\mathbb{C}[\hat{W}_{{M}'_\epsilon}\times
\hat{W}_{M''_\epsilon}]$ par restriction puis induction
tordue de fa\c con similaire le cas des groupes sym\'etriques;
c'est-\`a-dire que l'on fixe $M'_\epsilon,M''_\epsilon$ avec
$M'_\epsilon+M''_\epsilon=N'_\epsilon+N''_\epsilon$ et
consid\`ere les quadruplets $N_\epsilon^{i,j}$ o\`u $i,j\in
\{',''\}$ v\'erifiant
$N_\epsilon^{i,'}+N_\epsilon^{i,''}=N^i_\epsilon$ pour $i='$
et $''$ et
$N_\epsilon^{',j}+N_\epsilon^{'',j}=M^j_\epsilon$ pour $j='$
et $''$, s'il en existe (sinon on ne fait rien pour ce choix
de
$M_\epsilon',M_\epsilon''$). On restreint la repr\'esentation
de
$W_{\tilde{N'}_\epsilon}\times W_{\tilde{N}''_\epsilon}$ au
groupe $\times_{i,j}W_{N_\epsilon^{i,j}}$, on tensorise la
restriction par le caract\`ere de ce groupe qui vaut:
$sgn_{CD}^{(1-\zeta_\epsilon)/2}$ sur $W_{N_\epsilon^{',''}}$,
$\tilde{\chi}$ sur
$W_{N_\epsilon^{'','}}$,
$1$ sur $W_{N_\epsilon^{','}}$ et
$\tilde{\chi}sgn_{CD}^{(1+\zeta_\epsilon)/2}$ sur
$W_{N_\epsilon^{'',''}}$.

On induit au groupe $W_{M_\epsilon'}\times W_{M_\epsilon''}$
apr\`es avoir
\'echang\'e les 2e et 3e facteurs, c'est-\`a-dire
$W_{N_\epsilon^{',''}}$ et $W_{N_\epsilon^{'','}}$. Puis on
somme sur tous les quadruplets.  Ensuite on identifie
$\oplus_{(M_\epsilon',M_\epsilon'')\vert
M_\epsilon'+M_\epsilon''=N'_{\epsilon}+N''_\epsilon}
\mathbb{C}[\hat{W}_{{M}'_\epsilon}\times
\hat{W}_{M''_\epsilon}]$ \`a l'ensemble des combinaisons
lin\'eaires de base les couples de symboles ordonn\'es le
premier de d\'efaut impair
\'egal \`a $I_\epsilon$ et le deuxi\`eme de d\'efaut pair
\'egal \`a $\tilde{\zeta}_\epsilon P_\epsilon$.

C'est la construction de \cite{mw} 3.18 que l'on a
compl\`etement explicit\'ee. C'est assez compliqu\'e;
remarquons que la pr\'esence du caract\`ere $\tilde{\chi}$
n'a jou\'e de r\^ole dans \cite{mw} qu'en 5.5. Il en est de
m\^eme ici, ce caract\`ere ne joue aucun r\^ole sauf dans
l'\'enonc\'e de la conjecture \ref{conjecture}. La prise en
compte de $\zeta_\epsilon$ dans la d\'efinition, elle joue un
r\^ole mais dans la d\'efinition de $k_\chi$,
$\zeta_\epsilon$ aussi joue un r\^ole et en fait ces 2 prises
en compte se compensent en grande partie (cf. la preuve de
\ref{cle}).

\section{Localisation}

\subsection{Localisation des faisceaux
caract\`eres\label{localisation}}

On va avoir besoin d'une formule due \`a Lusztig qui calcule
les faisceaux caract\`eres au voisinage des points
semi-simples en terme de fonctions de Green. Elle est
\'ecrite en toute g\'en\'eralit\'e dans \cite{lusztig} et
explicit\'ee dans certains cas dans \cite{w3} et \cite{mw};
c'est la pr\'esentation de \cite{w3} par.7 que l'on reprend.
On fixe un \'el\'ement semi-simple $g_s$ de $SO(2n+1,F)\sharp$
et on suppose que toutes les valeurs propres de $g_s$ sont des
racines de l'unit\'e d'ordre premier \`a $p$. On fixe aussi
$g_u$ un \'el\'ement topologiquement unipotent de
$SO(2n+1,F)\sharp$ commutant \`a $g_s$. On suppose qu'il
existe un parahorique $K_{n',n''}$ (pour $n',n''$
convenables) contenant $g_sg_u$ et on note $s,u$ les
r\'eductions de $g_s$ et $g_u$ modulo le radical
pro-p-unipotent. 
On se donne aussi $\underline{n}\in D(\chi)$ que l'on suppose
relatif \`a $(n',n'')$ au sens 
$\sum_{[u]\in [VP(\chi)]}n'([u])\ell_{[u]}=n'$ et
$\sum_{[u]\in [VP(\chi)]}n''([u])\ell_{[u]}=n''$.

Dans la suite, on fixe, pour $\epsilon\in \{\pm 1\}$, des
donn\'ees $I_\epsilon,P_\epsilon, \tilde{\zeta}_\epsilon$
comme dans les paragraphes pr\'ec\'edents, c'est-\`a-dire une
donn\'ee cuspidale
$\tilde{cusp}$; on a choisi cette notation pour qu'elle soit
analogue \`a celle de la fin de \ref{inductionrestriction}.
Donc en particulier, la propri\'et\'e de
$\tilde{\zeta}_\epsilon$ est de v\'erifier
$\tilde{\zeta}_\epsilon=+$ si $P_\epsilon=0$. Et dans l'espace
vectoriel
$$
\oplus_{\epsilon\in \{\pm
1\}}\oplus_{m'(\epsilon),m''(\epsilon);
m'(\epsilon)+m''(\epsilon)=m(\epsilon)}\mathbb{C}[\hat{\cal
{W}}_{m'(\epsilon)}]\otimes \mathbb{C}[\hat{\cal
{W}}_{m'(\epsilon)}]$$ on ne regarde que les sous-espaces
vectoriels correspondant aux symboles relatifs
\`a ces donn\'es cuspidales. Pour manifester cette
restriction, on ajoute $\tilde{cusp}$ en indice. En fonction
de ce que l'on a rappel\'e en
\ref{bis} cela revient au m\^eme que de regarder les
repr\'esentations des diff\'erents groupes
$\times_{\epsilon\in \{\pm 1\}}W_{N'(\epsilon)}\times
W_{N''(\epsilon)}$, o\`u $W$ est un groupe de Weyl de type C
et o\`u les nombres
$N'(\epsilon),N''(\epsilon)$ v\'erifient les conditions:
$$ 2N'(\epsilon)+2N''(\epsilon)+(I_\epsilon^2+P_\epsilon^2)/2
=m(\epsilon).
$$On utilise la convention que si
$P_\epsilon=0$ alors
$\tilde{\zeta}=+1$. On pose
$\zeta_\epsilon=(-1)^{(I_\epsilon-1)/2}\tilde{\zeta}_\epsilon$
et on retrouve la convention de \ref{inductionrestriction}
que si
$P_\epsilon=0$, alors $
\zeta_\epsilon=(-1)^{(I_\epsilon-1)/2}$.

Soit alors $\phi$ un \'el\'ement de  ${\mathbb
C}[\hat{\cal{W}}_{\underline{n}, \tilde{cusp}}]$. Le point est de
calculer
$k_{\chi}(\phi)(g_sg_u)$ \`a l'aide des fonctions de Green du
commutant de $s$ dans le groupe $K_{n',n''}$ en r\'eduction.
Pour le faire, on suppose $s$ elliptique.

On d\'ecrit d'abord le commutant de
$g_s$ dans
$SO(2n+1,F)$; on note $[VP(g_s)]$ l'ensemble des valeurs
propres de $g_s$ regroup\'ees en paquets $\lambda$ et
$\lambda'$ sont dans le m\^eme paquet s'il existe $a\in 
\mathbb{N}$ tel que
$\lambda'=\lambda^{q^a}$. On note $m([\lambda])$ la
multiplicit\'e de $\lambda$. Si $\lambda \neq \pm 1$, on note
$\ell_{[\lambda]}:= 1/2\, \vert [\lambda]\vert$. Pour unifier
pour $\epsilon=\pm$, on pose $\ell_{\epsilon}=1$. Pour
$\lambda\neq \pm 1$, on note $F_{2\ell_{[\lambda]}}$
l'extension non ramifi\'ee de $F$ de degr\'e
$2\ell_{[\lambda]}$ et il existe une forme hermitienne (pour
l'extension $F_{2\ell_{[\lambda]}}/F_{\ell_{[\lambda]}}$)
$<\, ,\,>_{[\lambda]}$ sur l'espace vectoriel sur
$F_{2\ell_{[\lambda]}}$ de dimension $m([\lambda])$ tel que la
partie du commutant de $g_s$ relative \`a la valeur propre
$\lambda$ soit pr\'ecis\'ement le groupe unitaire de cette
forme. Des formes hermitiennes, comme ci-dessus, il y en a
exactement 2 qui se distinguent par la parit\'e de la
valuation du d\'eterminant, donc par un signe que nous
noterons
$\epsilon_{[\lambda]}$. Si l'on note
$U_{\epsilon_{[\lambda]}}$ le groupe de la forme
correspondant \`a $\epsilon_{[\lambda]}$, on rappelle que
$U_{\epsilon_{[\lambda]}}\simeq U_{-\epsilon_{[\lambda]}}$ si
$m([\lambda])$ est impair; mais nos constructions
d\'ependront de $\epsilon_{[\lambda]}$ comme en \cite{mw} 3.3
et suivants, et il faut donc garder la distinction.

Une valeur propre dans
$\{\pm 1\}$ introduit une forme orthogonale $<\, ,\, >_\pm$
($\pm$
 est ici le signe de la valeur propre consid\'er\'ee) sur un
$F$-espace vectoriel de dimension 
$m(\pm 1)$ et la partie du commutant qui lui correspond est
le groupe orthogonal de la forme; on utilise  la notation
$\eta_{\pm}$ et $\epsilon_{\pm}$ pour le signe du
discriminant et l'invariant de Hasse; il y a toujours un
probl\`eme sur la normalisation du discriminant et ici on
suit les conventions de \cite{mw}, c'est-\`a-dire que le
discriminant est invariant par ajout de plans hyperboliques
mais il n'est donc pas additif. Le signe du discriminant est
l'image du discriminant par le caract\`ere quadratique non
ramifi\'e de $F^*$ et la non additivit\'e n'est un probl\`eme
pour le signe du discriminant que si $-1$ n'est pas un
carr\'e.

Remarquons tout de suite, puisque l'on en aura besoin, que
toutes les quantit\'es qui viennent d'\^etre introduites sont
constantes sur la classe de conjugaison stable de
$g_s$ sauf la famille des $\epsilon_{[\lambda]}$ pour
$[\lambda]\in [VP(g_s)]$. Cette famille est   soumise \`a la
condition
$\times_{[\lambda]}\epsilon_{[\lambda]}=\sharp$ et l'ensemble
de ces familles soumises \`a cette condition param\'etrise
l'ensemble des classes de conjugaison dans la classe de
conjugaison stable de $g_s$ dans $SO(2n+1,F)_{\sharp}$. Si
l'on enl\`eve la condition de produit, l'ensemble plus grand
param\'etrise la classe de conjugaison stable de $g_s$ dans
$SO(2n+1,F)_{iso}\cup SO(2n+1,F)_{an}$.

Revenons au parahorique $K_{n',n''}$ que l'on a fix\'e au
d\'ebut et qui contient $g_s$ et $g_u$. On d\'efinit $s,u$
comme ci-dessus. D\'ecrivons le  commutant de $s$ dans la
partie r\'eductive du parahorique; on \'ecrit $s=(s',s'')$
avec $s'\in SO(2n'+1,\mathbb{F}_q)$ et $s''\in
O(2n'',\mathbb{F}_q)$. Les valeurs propres de $s$ sont les
''m\^emes`` que celles de $g_s$ ; pour $[\lambda]\in
[VP(g_s)]$, on note $m'([\lambda])$ la multiplicit\'e d'un
\'el\'ement
$\lambda\in [\lambda]$ comme valeur propre de $s'$ et
$m''([\lambda])$ l'analogue pour
$s''$. On a $m'([\lambda])+m''([\lambda])=m([\lambda])$. De
m\^eme pour $[\lambda]\in [VP(g_s)]$, $\lambda\neq \pm 1$, il
existe une forme hermitienne $<\, ,\, >_{[\lambda]}'$ (resp.
une forme hermitienne $<\, ,\, >_{[\lambda]}''$) sur un
$\mathbb{F}_{2\ell_{[\lambda]}}$ espace vectoriel 
de dimension
$m'([\lambda])$ (resp. $m''([\lambda])$) et des formes
orthogonales pour $\lambda=\pm 1$, telles que le commutant de
$s'$ (resp. $s''$) dans $SO(2n'+1,\mathbb{F}_q)$ (resp.
$O(2n'',\mathbb{F}_q)$) soit les \'el\'ements de
d\'eterminant 1 dans le produit (resp. le produit) des groupes
de ces formes. Il y a \'evidemment des rapports entre
$<\,,\,>'_{[\lambda]},<\,,\,>''_{[\lambda]}$ et $<\, ,\,
>_{[\lambda]}$. Supposons d'abord que  $\lambda\notin \{\pm
1\}$. Si $\epsilon_{[\lambda]}=1$ (d\'efini ci-dessus), alors
$m''([\lambda])$ est n\'ecessairement pair alors que ce
nombre est impair si
$\epsilon_{[\lambda]}=-1$, ensuite
$<\,,\,>'_{[\lambda]}\otimes <\,,\,>''_{[\lambda]}$ est
obtenu par r\'eduction (\`a l'aide d'un r\'eseau convenable)
de $<\, ,\, >_{[\lambda]}$.

Supposons maintenant que $\lambda\in \{\pm 1\}$. Alors
$m''([\lambda])$ a la m\^eme parit\'e que
$v_F(\eta_{[\lambda]})$. On note $\eta'_{[\lambda]}$ et
$\eta''_{[\lambda]}$ les discriminants des formes $<\, ,\,
>'_{[\lambda]}$ et $<\, ,\,>''_{[\lambda]}$; on les voit comme
des signes, c'est-\`a-dire qu'au lieu de regarder le
discriminant comme un \'el\'ement de $\mathbb{F}_q$ modulo
les carr\'es, on regarde sont image dans $\{\pm 1\}$. Si
$v_F(\eta_{[\lambda]})$ est pair, l'image de
$\eta''_{[\lambda]}$ dans $\mathbb{F}_q^*/\mathbb{F}_q^{*2}$
est
$\epsilon_{[\lambda]}$ tandis que si $v_F(\eta_{[\lambda]})$
est impair c'est l'image de $\eta'_{[\lambda]}$ dans le
m\^eme groupe qui est $\epsilon_{[\lambda]}$. De plus le
signe du discriminant  de $<\, ,\, >_{[\lambda]}$ v\'erifie:

$sgn(\eta_{[\lambda]})=sgn(-1)^{m'([\lambda])m''([\lambda])}
(-1)^
{v_F(\eta_{[\lambda]})}
\epsilon_{[\lambda]}\eta^{i}_{[\lambda]},$ o\`u
$i([\lambda])='$ si $v_F(\eta_{[\lambda]})$ est pair et
$''$ sinon.

On pose ici pour $[\lambda]\neq [\pm 1]$, ${\cal
W}_{\underline{m}_{g_s}(\lambda)}:=\mathfrak{S}_
{m'_{g_s}(\lambda)}\times \mathfrak{S}_{m''_{g_s}(\lambda)}$.
Pour $\lambda=\pm 1$, il y a des difficult\'es li\'ees \`a
l'existence de faisceaux cuspidaux; ici on ne s'int\'eresse
qu'aux fonctions de Green et les param\`etres pour la partie
cuspidales  sont alors des quadruplets d'entiers positifs ou
nuls. On \'ecrit les choses comme on en aura besoin; on avait
fix\'e ci-dessus une donn\'ee cuspidale,
$\tilde{cusp}$; on note $\vert cusp \vert$ un quadruplet
d'entiers positifs ou nuls,  
$\vert r\vert ^i_{\epsilon'}$ pour
$i\in
\{',''\}$ et $\epsilon'\in \{\pm 1\}$, qui est la partie
cuspidale pour les fonctions de Green.  Ce quadruplet
d\'epend de $\tilde{cusp}$ par les formules:

$ \vert r\vert '_{+1}, \vert r \vert '_{-1}$ est \`a l'ordre
pr\`es le couple
$I_++I_-, \vert I_+-I_-\vert$ avec $\vert r\vert'_{+1}$
impair par hypoth\`ese et $\vert r\vert ''_{+1}, \vert r\vert
''_{-1}$ est \`a l'ordre pr\`es le couple $P_++P_-, \vert
P_+-P_-\vert$ avec $\vert r\vert ''_{+1}\geq \vert r\vert
''_{-1}$ si et seulement si
$\zeta_+\zeta_-=(-1)^{1+(I_++I_-)/2}$.

On pose alors pour
$\lambda\in
\{\pm 1\}$ que l'on note plut\^ot $\epsilon'$, et pour un
couple d'entier $m'(\epsilon'),m''(\epsilon')$ v\'erifiant
$m'(\epsilon')+m''(\epsilon')=m(\epsilon')$
$$ {\cal W}_{m'(\epsilon'),m''(\epsilon'),\vert cusp \vert}:=
W_{1/2(m'(\epsilon')-\vert r'_{\epsilon'}\vert^2)}\times
W_{1/2(m''(\epsilon')-\vert r''_{\epsilon'}\vert^2)}$$\'etant
entendu que ce groupe est nul si l'un des indices n'est pas
un entier positif ou nul.

Pour la donn\'ee d'un ensemble de couple
$\underline{m}_{g_s}:=\{m'(\lambda),m''(\lambda)\in
D(m_{g_s}(\lambda))\}$, on pose ${\cal
{W}}_{\underline{m}_{g_s}, \vert cusp\vert}$ le produit des
groupes d\'efinis ci-dessus. Et les fonctions de Green donnent
une application de
${\mathbb C}[\hat{{\cal{W}}}_{\underline{m}_{g_s}}]$ dans
l'ensemble des fonctions \`a support unipotent du commutant
de $g_s$ dans $SO(2n'+1,\mathbb{F}_q)\times
O(2n'',\mathbb{F}_q)$; cela s'interpr\`ete encore comme des
fonctions
\`a support dans l'ensemble des \'el\'ements topologiquement
unipotents de $K_{n',n''}$ commutant \`a
$g_s$. Quand on somme sur tous les $\underline{m}_{g_s}$, on
d\'efinit $\mathbb{C} [\hat{\cal W}_{D(g_s),\vert cusp\vert}]$
et on s'autorisera la suppression du $cusp$ quand il n'y a pas
d'ambiguit\'e. 

On revient maintenant \`a $\phi$ et donc pour tout $[u]\in
[VP(\chi )]$ on a un
\'el\'ement $(n'[u],n''[u])\in D(m([u]))$ de telle sorte que
$\sum_{[u]\in [VP(\chi)]} n'[u]\ell_{[u]}=n'$ et on rappelle
que l'on a aussi fix\'e une donn\'ee cuspidale,
$\tilde{cusp}$. On note
$\underline{n}$ l'ensemble de ces paires et on a d\'efini
${\cal{W}}_{\underline{n}, \tilde{cusp}}$. On va d\'efinir une
application de ${\mathbb C}[\hat{{\cal{W}}}_{\underline{n}}]$
dans
${\mathbb C}[\hat{{\cal{W}}}_{\underline{m}_{g_s}}] $. Mais
on a besoin d'un certain nombre d'objets interm\'ediaires.

Soit une collection de paires
$\underline{\nu}([u],[\lambda])=
(\nu'([u],[\lambda]),\nu''([u],[\lambda]))$ soumises aux
conditions, o\`u $\ell_{[u]}$ et $\ell_{[\lambda]}$ sont
comme ci-dessus (pour $a,b$ des entiers, on note $(a,b)$ le
pgcd de ces nombres):
$$\sum_{[u]}\tilde{\nu}
([u],[\lambda])\ell_{[u]}/(\ell_{[u]},\ell_{[\lambda]}) =
\underline{m}_{g_s}([\lambda]);$$ la somme de couples,  se
fait terme \`a terme. Et on a aussi:
$$\sum_{[\lambda]}\nu([u],[\lambda])
/(\ell_{[u]},\ell_{[\lambda]})=
\underline{n}([u]).$$ On pose
${\cal{W}}_{\nu([u],[\lambda])}:=
\mathfrak{S}_{\nu'([u],[\lambda])} \times
\mathfrak{S}_{\nu''([u],[\lambda])} $ si soit $[u]$ soit
$[\lambda]$ n'est un \'el\'ement de $\{\pm 1\}$. Pour
traiter le cas de $\pm 1$, on r\'eutilise la donn\'ee de la
partie cuspidale $\vert cusp\vert$ et on pose, pour
$\epsilon\in
\{\pm 1\}$ et $\epsilon'\in  \{\pm 1\}$
${\cal{W}}_{\underline{\nu}([\epsilon],[\epsilon'])}=$ 

$W_{1/2(\nu'([\epsilon],[\epsilon'])-(\vert
r\vert '_{\epsilon'})^2)}
\times
{W}_{1/2(\nu''([\epsilon],[\epsilon'])-(\vert
r\vert ''_{\epsilon'})^2)}$; en particulier cela sous-entend
que ces nombres \'ecrits en indice sont des entiers positifs
ou nuls (sinon le groupe d\'efini est 0). Et on note
${\cal{W}}_{\underline{\nu}(\chi,g_s),cusp}$ l'union de tous
ces groupes.

Pour $\delta='$ ou $''$ et pour $[u]\in [VP(\chi)]$,
$[\lambda]\in [VP(g_s)]$, on a besoin de d\'efinir une
application de ${\cal W}_{\nu^\delta([u],[\lambda])}$ dans
${\cal W}_{n^\delta([u])_{[\lambda]}}$ et dans ${\cal
W}_{n^\delta([\lambda])_{[u]}}$ on pr\'ecisera dans chaque
cas les rapports entre les entiers
$\nu^\delta([u],[\lambda]), n^\delta([u])_{[\lambda]}$ et
$\nu^\delta([\lambda])_{[u]}$. Cette application n'est pas un
morphisme de groupes mais simplement compatible
\`a l'action adjointe; au passage, on d\'efinit aussi une
fonction invariante  par conjugaison sur ${\cal
W}_{\nu^\delta([u],[\lambda])}$ que l'on note
$\chi^\delta_{[u],[\lambda]}$. Le point qui n'est pas nouveau
est que pour les groupes unitaires qui interviennent soit
pour $[u]$ quand $[u]\neq \pm 1$ soit pour $[\lambda]$ quand
$[\lambda]\neq \pm 1$, un tore n'est pas vraiment associ\'e
\`a un \'el\'ement du groupe sym\'etrique convenable mais au
produit d'un tel \'el\'ement par un Frob\'enius.

On d\'ecrit ces constructions au cas par cas: premier cas:
$[u], [\lambda] \in \pm 1$; ici, on veut
$\nu^\delta([u],[\lambda])=n^\delta([u])_{[\lambda]}$ et 
l'application de
$ {\cal W}_{\nu^\delta([u],[\lambda])}={
W}_{\nu^\delta([u],[\lambda])}$ dans ${\cal
W}_{n^\delta([u])_{[\lambda]}}={
W}_{n^\delta([u])_{[\lambda]}}$ est l'identit\'e. On fait une
construction analogue en rempla\c cant $[u]$ par $[\lambda]$.
Quand
\`a la fonction
$\chi^\delta_{[u],[\lambda]}$ c'est l'identit\'e sauf si
$u=\lambda=-1$ o\`u c'est le
$sgn_{CD}(-1)^{\nu^\delta([u],[\lambda])(q-1)/2}$. 

deuxi\`eme cas: $[u]\neq [\pm 1]$, $[\lambda] \neq [\pm 1]$
et $\ell_{[u]}$ et
$\ell_{[\lambda]}$ sont divisibles par la m\^eme puissance de
2 (cette derni\`ere condition donne des renseignements sur
les tores des groupes unitaires intervenant). L'application de
${\cal W}_{\nu^\delta([u],[\lambda])}= {\mathfrak
S}_{\nu^\delta([u],[\lambda])}$ dans ${\mathfrak
S}_{n^\delta([u])_{[\lambda]}}$ se d\'ecrit quand
$n^\delta([u])_{[\lambda]}=\nu^\delta([u],[\lambda])\ell_{[\lambda]}/(\ell_{
[u]},\ell_{[\lambda]})$. A une permutation dont les cycles
sont $(\alpha_1,\cdots, \alpha_r)$ on associe la permutation
de cycles $(\alpha_1\ell_{[\lambda]}/(\ell_{
[u]},\ell_{[\lambda]}), \cdots, \alpha_r
\ell_{[\lambda]}/(\ell_{ [u]},\ell_{[\lambda]}))$. Quand on
consid\`ere $[\lambda]$ au lieu de $[u]$, on \'echange
simplement les r\^oles de
$\ell_{[u]}$ et $\ell_{[\lambda]}$. Pour d\'ecrire  la
fonction
$\chi^\delta_{[u],[\lambda]}$, on a besoin d'une notation
auxiliaire. On pose, pour $\alpha \in \mathbb{N}$, $y\in
\mathbb{Q}$:
$$\phi_{\alpha,[\lambda],y}:=y^{-1}(\ell_{[\lambda]})^{-1}(
\lambda^{\alpha }+\lambda^{\alpha q}+\cdots +\lambda^{\alpha
q^{\ell_{[\lambda]}-1}}+\lambda^{-\alpha}+\lambda^{-\alpha
q}+\cdots +\lambda^{-\alpha q^{\ell_{[\lambda]}-1}}),$$o\`u
$\lambda$ est n'importe quel \'el\'ement dans $[\lambda]$.
Alors  sur l'\'el\'ement $w$ associ\'e aux cycles $\alpha_1,
\cdots ,\alpha_r$,
$\chi^\delta_{[u],[\lambda]}(w)=\prod_{s\in
[1,r]}\phi_{\alpha_s,[\lambda],(\ell_{[u]},\ell_{[\lambda]})}$.

troisi\`eme cas: $[u]\neq [\pm 1]$, $[\lambda] \neq [\pm 1]$
et $\ell_{[u]}$ et
$\ell_{[\lambda]}$ ne sont pas divisibles par la m\^eme
puissance de 2. Ici on veut alors,
$n^\delta([u])_{[\lambda]}=2\nu^\delta([u],[\lambda])\ell_
{[\lambda]}/(\ell_{[u]},\ell_{[\lambda]})$. A une permutation
dont les cycles sont $(\alpha_1,\cdots, \alpha_r)$ on
associe, ici, la permutation de cycles
$(2\alpha_1\ell_{[\lambda]}/(\ell_{ [u]},\ell_{[\lambda]}),
\cdots, 2\alpha_r
\ell_{[\lambda]}/(\ell_{ [u]},\ell_{[\lambda]}))$. Quand on
consid\`ere $[\lambda]$ au lieu de $[u]$, on \'echange
simplement les r\^oles de
$\ell_{[u]}$ et $\ell_{[\lambda]}$. Et la fonction
$\chi^\delta([u],[\lambda])$ vaut sur ce $w$, $
\prod_{s\in [1,r]}\phi_{\alpha_s,[\lambda],
(\ell_{[u]},\ell_{[\lambda]})/2}$.

quatri\`eme cas: $u=\pm 1$ et $[\lambda]\neq [\pm 1]$. Ici on
veut
$n^\delta([u])_{[\lambda]}=\nu^\delta([u],[\lambda])
\ell_{[\lambda]}$ et
$n^\delta([\lambda])_{[u]}=\nu^\delta([u],[\lambda])$.
L'application de $\mathfrak{S}_{\nu^ \delta([u],[\lambda])}$
dans $\mathfrak{S}_{n^\delta([\lambda])_{[u]}}$ est
l'identit\'e. L'application de 
$\mathfrak{S}_{\nu^ \delta([u],[\lambda])}$ dans
$W_{n^\delta([u])_{[\lambda]}}$ envoie l'\'el\'ement $w$ de
cycle $(\alpha_1,\cdots,
\alpha_{r},\alpha'_1,\cdots, \alpha'_{r'})$, o\`u les
$\alpha$ sont pairs et les $\alpha'$ impairs sur l'\'el\'ement
de $W_{n^\delta([u])_{[\lambda]}}$ qui correspond \`a 2
partitions $(\alpha_1 \ell_{[\lambda]}, \cdots,
\alpha_r\ell_{[\lambda]}; \alpha'_1 \ell_{[\lambda]}, \cdots,
\alpha'_{r'}\ell_{[\lambda]})$. Quand \`a la fonction
$\chi^\delta([u],[\lambda])$, c'est l'identit\'e si
$[u]=[+1]$ et est constante de valeur
$(-1)^{\nu^\delta([u],[\lambda])(1+q^{\ell_{[\lambda]}})/2}$
pour
$[u]=[-1]$.

cinqui\`eme cas: $[u]\neq [\pm 1]$ et $\lambda=\pm 1$. On
\'echange les r\^oles de $[u]$ et $[\lambda]$ dans le cas
ci-dessus.

Il faut aussi prendre en compte une contribution de la partie
cuspidale; l\`a il n'y a pas de groupes mais simplement une
fonction \`a d\'efinir, $\chi^\delta_{cusp,\underline{\nu}}$
(o\`u $\delta\in \{',''\}$ comme ci-dessus) qui d\'epend de
$\underline{\nu}\in D(\chi,g_s)$ et de $w\in {\cal
W}_{\underline{\nu}}$. Pour cela, on pose
$\chi^\delta_{+,\underline{\nu}}(w):=(-1)^{\sum_{[u]\in
[VP(\chi)]-\{\pm 1\}}\nu^\delta([u],+)}
sgn_{CD}w_{\nu^\delta(+,+)}w_{\nu^\delta(-,+)}$ et une
d\'efinition analogue, $\chi^\delta_{-,\nu}$ en rempla\c cant
$+$ par
$-$.  On note aussi
$\epsilon_I$ (resp. $\epsilon_P$) le signe tel que
$I_{\epsilon_I} -I_{-\epsilon_I}>0$ (resp. $P_{\epsilon_P} -
P_{-\epsilon_P}>0$); si les 2 nombres sont
\'egaux \`a 0, on prend le signe $\epsilon_I=\epsilon_P=+$ par
convention et si l'un seulement des nombres vaut 0, alors on
prend par convention $\epsilon_I=\epsilon_P$. On pose
$r'_{im}$ l'\'el\'ement impair du couple
$(I_++I_-)/2, \vert I_+-I_-\vert/2$ et l'on a ($x'\in
\mathbb{N}$ ne nous sert \`a rien mais est le
$\delta(r',r'')$ de \cite{w3} et il en est de m\^eme pour
$x''$):
\[\chi'_{cusp,\underline{\nu}}(w)=q^{x'}(-1)^{(q-1)(r'-1)/4}
(\chi'_{+,\underline{\nu}}(w)
\chi'_{-,\underline{\nu}}(w)\eta'_+(g_s)\eta'_-(g_s))^
{(I_{\epsilon_I}-1)/2}
\begin{cases}1 \hbox{ si }\epsilon_I=+\\
\eta'_-(g_s)\chi'_{-,\underline{\nu}}(w) \hbox{ si
}\epsilon_I=-.\end{cases}
\] ci-dessous $y''$ est  $r''_+$ si ce nombre est pair,
et
$r''_--1$ sinon,
\[\chi''_{cusp,\underline{\nu}}(w)=q^{x''}(-1)^(y''(q-1)/4)
(\eta''_+(g_s)\eta''_-(g_s)\chi''_{+,
\underline{\nu}}\chi''_{-,
\underline{\nu}})^{(I_{\epsilon_P}-1)/2}
\]
\[
\begin{cases}1\hbox{ si }{\zeta}_{\epsilon_P}>0 \hbox{ et
}\epsilon_P=+\\
\eta''_-(g_s)\chi''_{-,\underline{\nu}}(w) \hbox{ si
}{\zeta}_{\epsilon_P}=+ \hbox{ et }\epsilon_P=-\\
\eta''_+(g_s)\chi''_{+,\underline{\nu}}(w) \hbox{ si
}{\zeta}_{\epsilon_P}=- \hbox{ et } \epsilon_P=-\\
\eta''_+(g_s)\eta''_-(g_s)\chi''_{+,\underline{\nu}}(w)
\chi''_{-,\underline{\nu}}(w)\hbox{ si }{\zeta}_{\epsilon_P}=-
\hbox{ et }\epsilon_P=+.
\end{cases}
\] Le produit de ces 2 fonctions se simplifie un peu. On
\'ecrit ce produit comme le produit des 3 termes

$C_{cusp}$ qui est une constante,

$c_{cusp}(g_s):=
(\eta'_+(g_s)\eta'_-(g_s))^{(I_{\epsilon_I}-1)/2}
(\eta''_+(g_s)\eta''_-(g_s))^{(I_{\epsilon_P}-1)/2} 
\times$
$$
\begin{cases} 1\text{ si } \zeta_{\epsilon_P}=+\\
\eta'_+(g_s)\eta'_-(g_s)\text{ si }\zeta_{\epsilon_P}=-
\end{cases}
\times
\begin{cases}1\text{ si }\epsilon_I=+\\\eta'_-(g_s)\text{ si
}\epsilon_I=-
\end{cases}
\times
\begin{cases} 1\text{ si }\epsilon_P=+\\\eta''_-(g_s)\text{ si
}\epsilon_P=-
\end{cases}$$

$
\chi_{cusp,\underline{\nu}}(w):= (\chi'_{+,\underline{\nu}}(w)
\chi'_{-,\underline{\nu}}(w))^{(I_{\epsilon_I}-1)/2}
(\chi''_{+,\underline{\nu}}(w)
\chi''_{-,\underline{\nu}}(w) )^{(I_{\epsilon_P}-1)/2} 
\times$
$$
\begin{cases} 1\text{ si } \zeta_{\epsilon_P}=+\\
\chi'_{+\underline{\nu}}(w)
\chi'_{-,\underline{\nu}}(w)\text{ si }\zeta_{\epsilon_P}=-
\end{cases}
\times
\begin{cases}1\text{ si }\epsilon_I=+\\
\chi'_{-,\underline{\nu}}(w)
 \text{ si }\epsilon_I=-
\end{cases}
\times
\begin{cases} 1\text{ si }\epsilon_P=+\\
\chi''_{-,\underline{\nu}}(w)\text{ si }\epsilon_P=-
\end{cases}$$

Supposons maintenant donn\'e $\underline{n}\in D(\chi)$ et
$\underline{m}\in D(g_s)$. Soit $\underline{\nu}\in
D(\chi,g_s)$ et on suppose que  les
\'egalit\'es suivantes sont v\'erifi\'ees:

pour
$\delta\in
\{',''\}$, pour tout $[u]\in [VP(\chi)]$, 
$\sum_{[\lambda]\in [VP(g_s)]}
\nu^\delta([u],[\lambda])2^{x_{[u],[\lambda]}}\ell_{[\lambda]}
/(\ell_{[u]},\ell_{[\lambda]})=n^\delta([u]),$ o\`u
$x_{[u],[\lambda]}=0$ sauf si 
$\ell_{[u]}/(\ell_{[u]},\ell_{[\lambda]})$ ou
$\ell_{[\lambda]}/(\ell_{[u]},\ell_{[\lambda]})$ est pair,
o\`u il vaut 1.

Et pour $\delta\in \{',''\}$, pour tout $[\lambda]\in
[VP(g_s)]$,
$\sum_{[u]\in [VP(\chi)]}
\nu^\delta([u],[\lambda])2^{x_{[u],[\lambda]}}\ell_{[u]}
/(\ell_{[u]},\ell_{[\lambda]})=m^\delta([\lambda]),$ o\`u
$x_{[u],[\lambda]}$ est comme ci-dessus.

En faisant des produits convenables des constructions
ci-dessus, on a une application de ${\cal
W}_{\underline{\nu}}$ d'une part dans ${\cal
W}_{\underline{n}}$ et d'autre part dans ${\cal W}_{
\underline{m}}$. Partant donc d'un \'el\'ement de ${\mathbb
C}[\hat{\cal {W}}_{\underline{n}}]$ on peut le restreindre en
un \'el\'ement de ${\mathbb C}[\hat{\cal
{W}}_{\underline{\nu}}]$, le multiplier par le
$\chi'_{cusp,\underline{\nu}}\chi''_{cusp,\underline{\nu}}\prod_
{[u]\in [VP(\chi)],[\lambda]\in [VP(g_s)]}
\chi'([u],[\lambda])\chi''([u],[\lambda])$
 puis l'induire en un \'el\'ement de 
${\mathbb C}[\hat{\cal {W}}_{\underline{m}}]$.

 L'application cherch\'ee est  la somme sur toutes les
collections
$\underline{\nu}(\chi,g_s)$ comme ci-dessus. Elle est
not\'ee, $loc_{g_s,\underline{m}}$.

\bf D\'efinition: \sl pour $g_s$ comme ci-dessus, on note
$loc_{g_s}:=\sum_{\underline{m}\in
D(g_s)}loc_{g_s;\underline{m}}$.\rm

 On rappelle que
les fonctions de Green d\'efinissent, pour $\underline{m}\in
D(g_s)$,
 une application de
$\hat{\cal W}_{\underline{m}}$ dans les fonctions sur le
centralisateur de $s$ (la r\'eduction de
$g_s$) dans la r\'eduction de
$K_{n',n''}$; on remarque que $n'$ et $n''$ sont
d\'etermin\'es par $\underline{m}$. On remonte ensuite en des
fonctions sur un sous-groupe compact ouvert convenable du
centralisateur de
$g_s$ dans
$SO(2n+1,F)_{\sharp}$ par invariance sous le radical
pro-p-unipotent puis on les prolonge par z\'ero en des
fonctions sur ce centralisateur. On note
$Q$ cette application. On peut donc d\'efinir
$Q(loc_{g_s}(\rho))$  pour $\rho\in \hat{{\cal
W}}_{\underline{n}}$ o\`u $\underline{n} \in D(\chi)$ (avec
une donn\'ee cuspidale fix\'ee)

\bf Remarque: \sl la d\'efinition pr\'ec\'edente ne d\'epend
de $g_s$ dans sa classe de conjugaison stable que par le
facteur $c_{cusp}(g_s)$.

\rm Pour cette section, plut\^ot que de travailler avec
$k_{\chi}$ qui a \'et\'e normalis\'e  pour avoir de bonnes
propri\'et\'es (cf. \cite{w3}), on \'enonce un r\'esultat
pour un analogue de $k_{\chi}$ non normalis\'e,
c'est-\`a-dire que les fonctions de Green g\'en\'eralis\'ees
associent, pour $\underline{n}$ fix\'e dans $D(\chi)$ et une
donn\'ee cuspidale fix\'ee $\tilde{cusp}$, \`a un \'el\'ement
du groupe $w\in {\cal{W}}_{\underline{n},\tilde{cusp}}$ une
fonction sur un groupe fini convenable. Et, pour $\phi$ une
repr\'esentation de
${\cal{W}}_{\underline{n},\tilde{cusp}}$, on a d\'efini
$k_{\chi}(\phi):= \sum_{w\in
{\cal{W}}_{\underline{n},\tilde{cusp}}}\lambda(w)tr(\phi)(w)
k_{\chi}(w)$, o\`u $\lambda(w)$ est un caract\`ere qui
d\'epend de la donn\'ee cuspidale. On utilise ici simplement
$k_\chi^{nn}(\phi)$ l'analogue en supprimant $\lambda$. La
raison est que $\rho\circ \iota$ d\'epend aussi du support
cuspidal et que l'on ne s'int\'eresse qu'au compos\'e $k_\chi
\circ \rho\circ \iota$; les normalisations s'annulent
partiellement et il vaut donc mieux ne pas se fatiguer \`a
les faire.
Et on a, pour
$g_u$ un
\'el\'ement topologiquement unipotent qui commute \`a $g_s$:

\bf Lemme: \sl Pour $\tilde{cusp}$, $\underline{n}$ et $\phi$
comme ci-dessus, (i)$k^{nn}_{\chi}(
{\phi})(g_sg_u)= Q(loc_{g_s;\underline{m}}(\phi))(g_u) .$

(ii)On a les \'egalit\'es d'int\'egrales orbitales (o\`u le
groupe est mis en exposant)
$$I^{SO(2n+1,F)_{\sharp}}(g_s,k^{nn}_{\chi}
(\phi))=  I^{Cent^0_{SO(2n+1,F)}g_s}(g_u,Q(loc_{g_s}(
\phi)).$$
\rm (i) C'est un probl\`eme sur le groupe fini
$SO(2n'+1,\mathbb{F}_q)\times O(2n'',
\mathbb{F}_q)$ o\`u il s'agit de localiser au voisinage de la
r\'eduction de $g_s$ le faisceau caract\`ere associ\'e \`a
$\underline{\rho}$. Cela a \'et\'e fait en toute
g\'en\'eralit\'e par Lusztig et il faut expliciter ses
formules. En \cite{mw} 2.16, on a trait\'e le cas o\`u
$[VP(\chi)]$ est r\'eduit \`a $+1$ mais il n'y a pas de
restriction sur $g_s$. En
\cite{w3} est trait\'e le cas o\`u $[VP(\chi)]$ et $[VP(g_s)]$
contiennent
$+1$ et
$-1$. La formule de Lusztig s'applique aux faisceaux
caract\`eres associ\'es \`a des \'el\'ements de ${\cal
W}_{D(\chi)}$ et non pas aux repr\'esentations de ce groupe.
La d\'emonstration de
\cite{w3} 7.1 qui se place dans ce cadre est tr\`es
g\'en\'erale et elle montre que le seul point est le calcul
de la constante not\'ee $z_2$ en loc.cit. Cette constante est
la somme  des valeurs du caract\`ere d\'etermin\'e par $\chi$
sur les conjugu\'es de $s$. Le calcul est plus compliqu\'e
qu'en loc.cite  mais c'est le calcul des
$\chi'([u],[\lambda])$ et $\chi''([u],[\lambda])$. Ensuite
\cite{w3} 7.2  d\'eduit le r\'esultat cherch\'e.

(ii) pour pouvoir utiliser (i), on d\'ecompose l'orbite de
$g_s$ sous $SO(2n+1,F)_{\sharp}$ en orbites sous $K_{n',n''}$
(le parahorique qui sert \`a la d\'efinition de
$k_{\chi}(\underline{\rho})$). Ces orbites sont
param\'etr\'ees par les \'el\'ements de $D(g_s)$; seuls
comptent les orbites qui coupent $K_{n',n''}$ et pour cela il
faut la relation:
$$\sum_{[\lambda]\in
[VP(\chi)]}m'([\lambda])\ell_{[\lambda]}=n'.$$ Si cette
relation n'est pas satisfaite, on remarque que les sommes
intervenant sont vides et
$loc_{g_s;\underline{m}} (\underline{\rho})=0$. On n'a donc
pas \`a se pr\'eoccuper de cette condition. Ensuite on
calcule l'int\'egrale sous chaque $K_{n',n''}$-orbite en
utilisant (i). Il y a clairement des mesures \`a prendre en
compte; c'est fait en \cite{mw} 3.17, le $\vert W(d)\vert$
n'intervient pas pour nous car il a \'et\'e pris en compte
quand on travaille avec des repr\'esentations des groupes
${\cal W}$ et non les \'el\'ements de ces groupes et les
constantes ($c(\gamma),c(\gamma)_{\sharp}$ de loc.cit) ont
\'et\'e mises dans la d\'efinition de $loc_{g_s}$.

\subsection{Restriction des repr\'esentations et localisation
des faisceaux caract\`eres\label{cle}} 

Dans cette section, il s'agit de montrer que l'op\'eration de
restriction aux parahoriques des repr\'esentations commute
\`a l'action de restriction des caract\`eres aupr\`es des
\'el\'ements semi-simples elliptiques, m\^eme si l'on ne peut
l'exprimer en ces termes tant que la conjecture
\ref{conjecture} n'est pas d\'emontr\'ee. En plus comme on
peut s'y attendre, vu la complexit\'e des formules, il n'y a
vraiment commutation que dans les cas favorables.

On fixe une donn\'ee cuspidale, $cusp$, pour les faisceaux
caract\`eres, c'est-\`a-dire pour nous, 2 entiers impairs,
$I_+$ et $I_-$ et 2 entiers pairs $P_+,P_-$ ainsi que 2 signes
$\zeta_+,\zeta_-$ avec la convention que si pour
$\epsilon=\pm$, $P_\epsilon=0$ alors
$\zeta_\epsilon=(-1)^{(I_\epsilon-1)/2}$. Dans
${\mathbb C}[\hat{\cal{W}}_{D(\chi)}]$, on ne consid\`ere que
la partie relative \`a cette donn\'ee cuspidale; on d\'efinit
comme en \ref{inductionrestriction}  une autre donn\'ee
cuspidale
$\tilde{cusp}$ simplement en changeant les signes,
$\tilde{\zeta}_{\pm}:=(-1)^{(I_{\pm}-1)/2}\zeta_\pm$ et on
reprend la notation $\tilde{\chi}$ de loc.cite; on note
$\tilde{\empty}$ l'application \'evidente de la partie de
${\mathbb C}[\hat{\cal{W}}_{D(\chi)}]$ relative \`a $cusp$
dans son homologue relative \`a $\tilde{cusp}$ qui en terme
de repr\'esentation de groupes de Weyl est tout simplement
l'identit\'e (mais on a chang\'e le support cuspidal)
compos\'e avec la multiplication par le caract\`ere
$\tilde{\chi}$.

Soit $g_s$ un \'el\'ement semi-simple elliptique dont les
valeurs propres sont des racines de l'unit\'e d'ordre premier
\`a $p$. On reprend la notation
$[VP(g_s)]$ pour signifier l'ensemble des valeurs propres de
$g_s$ regroup\'ees en paquets $\lambda,\lambda'$ sont dans le
m\^eme paquet s'il existe $a\in \mathbb{N}$ tel que
$\lambda'=\lambda^{q^a}$; la multiplicit\'e d'une valeur
propre $\lambda$ est not\'ee $m([\lambda])$ car elle ne
d\'epend que du paquet auquel $\lambda$ appartient.

On note $D(g_s)$ l'ensemble des d\'ecompositions
$\{(m'([\lambda]),m''([\lambda])\in
D(m([\lambda]))\}_{[\lambda]\in [VP(g_s)]}$; pour chaque
\'el\'ement $\underline{m}\in D(g_s)$, on a une classe
d'association de parahorique $K_{n',n''}$ o\`u
$n'=\sum_{[\lambda]}m'([\lambda])$ et \`a l'int\'erieur de ce
parahorique une classe de conjugaison d'\'el\'ements
semi-simples de r\'eduction semi-simple elliptique incluse
dans la classe de conjugaison de $g_s$; on connait les
valeurs propres des \'el\'ements dans cette classe de
conjugaison. Pour
$\underline{m}\in D(g_s)$, on reprend la notation
${\cal{W}}_{\underline{m}}$ de \ref{localisation}. On a
d\'efini ci-dessus des op\'erations de localisation de
${\mathbb C}[\hat{{\cal{W}}}_{D(\chi)}]$ dans
${\mathbb C}[\hat{{\cal{W}}}_{\underline{m}}]$; en sommant
sur tous les \'el\'ements de $D(g_s)$, on d\'efinit donc une
application de localisation de ${\mathbb C}
[\hat{{\cal{W}}}_{D(\chi)}]$ dans
${\mathbb C}[\hat{{\cal{W}}}_{D(g_s)}]$ que l'on note
$loc_{g_s}$. 

On remarque que $\rho\circ \iota$ se d\'efinit aussi de
${\mathbb C}[\hat{{\cal{W}}}_{D(g_s)}]$ dans lui-m\^eme; ce
sont exactement les d\'efinitions de \cite{mw} 3.1, 3.2 et
3.9, 3.10. On les pr\'esente un peu diff\'eremment de fa\c
con  similaire aux formules de
\ref{inductionrestriction}. Pr\'ecis\'ement, on note
$DD(g_s)$ l'ensemble des u-plets
$\{m^{i,j}([\lambda]); i,j\in \{',''\}, [\lambda]\in
[VP(g_s)]$. Pour $\underline{\underline{m}}\in DD(g_s)$, on
d\'efinit $\hat{{\cal W}}_{\underline{\underline{m}}}$ de fa\c
con similaire \`a \ref{inductionrestriction}, la partie
cuspidale possible
\'etant simplement donn\'ee comme en \ref{localisation}
(not\'ee $\vert cusp\vert$). Et ici,
$\rho\circ \iota$ est la somme sur tous les
$\underline{\underline{m}}\in DD(g_s)$ de la restriction \`a
${\cal W}_{\underline{\underline{m}}}$ suivie de l'induction
apr\`es avoir tordu par $(-1)^{\sum_{[\lambda]\notin \{\pm
1\}} m^{'',''}([\lambda])}(sgn_{CD})_{\vert
W_{N_+^{'',''}}\times W_{N_-^{'',''}}}$, les notations et les
inclusions entre les groupes
\'etant celles d\'ecrites en
\ref{inductionrestriction}.

 Pour traiter tous les cas, on pose encore quelques
d\'efinitions. Pour  $\epsilon'\in \{\pm 1\}$, on note
$X_{\epsilon'}$ l'endomorphisme de
${\mathbb C}[\hat{\cal{W}}_{D(g_s)}]$ qui est la
tensorisation par le caract\`ere trivial sur tous les
facteurs sauf $\hat{\cal {W}}_{m''(\epsilon')}$ o\`u il vaut
$sgn_{CD}$, suivie par l'inversion des facteurs relatifs
$m'(\epsilon')$ et $m''(\epsilon')$ (\`a ce stade cette
inversion est assez formelle mais elle a de l'importance
quand ensuite on applique $\rho\circ \iota$). On rappelle la
donn\'ee cuspidale fix\'ee et on reprend les notations
$\epsilon_I$ et $\epsilon_P$ de \ref{localisation}. On note
$X_{cusp}$ l'endormophisme de ${\mathbb
C}[\hat{\cal{W}}_{D(g_s)}]$ d\'efini par:
\[X_{cusp}:=\begin{cases} 1 \text{ si } \epsilon_I=\epsilon_P
\text{ et }\zeta_{\epsilon_P}=+\\ X_+X_- \text{ si }
\epsilon_I=\epsilon_P
\text{ et }\zeta_{\epsilon_P}=-\\
X_{(-1)^{1+(I_{+}+I_{-})/2} }\text{ si }
\epsilon_I\neq\epsilon_P \text { et } \zeta_{\epsilon_P}=-\\
X_{(-1)^{(I_{+}+I_{-})/2} }\text{ si }
\epsilon_I\neq\epsilon_P \text { et } \zeta_{\epsilon_P}=+ .
\end{cases}
\]

\bf Lemme: \sl Fixons la donn\'ee cuspidale comme ci-dessus.
Alors le diagramme ci-dessous est commutatif pour tout
\'el\'ement
$g_s$ comme en
\ref{localisation}:
\[\begin{array}{lcc} {\mathbb C}[\hat{{\cal{W}}}_{D(\chi)}]
&\stackrel{
\rho\circ
\iota} {\rightarrow }&{\mathbb C}[\hat{{\cal{W}}}_{D(\chi)}]\\
\Big\downarrow  X_{cusp}\circ loc_{g_s}\circ \tilde{\empty}&
&\Big\downarrow loc_{g_s}\\ {\mathbb
C}[\hat{{\cal{W}}}_{D(g_s)}] &\stackrel{
\rho\circ
\iota} {\rightarrow } &{\mathbb C}[\hat{{\cal{W}}}_{D(g_s)}]
\end{array}
\]
\rm

Pour d\'emontrer ce lemme, on r\'eintroduit le groupe
auxiliaire
${\cal{W}}_{DD(\chi)}$ et son avatar
${\cal{W}}_{DD(g_s)}$ qui permettent de calculer
$\rho\circ \iota$. De m\^eme, on r\'eintroduit
${\cal{W}}_{D(\chi,g_s)}$. On esp\`ere que le lecteur voit
une  localisation de
${\cal{W}}_{DD(\chi)}$ vers ${\cal{W}}_{DD(g_s)}$ qui utilise
le groupe
${\cal{W}}_{DD(\chi,g_s))}$ sugg\'er\'ee par les notations;
ce groupe est construit comme tous les groupes de m\^eme type
mais en utilisant des collections d'entiers
$\nu^{i,j}([u],[\lambda])$ pour $i,j\in
\{',''\}$, $[u] \in [VP(\chi)]$ et $[\lambda]\in [VP(g_s)]$
qui v\'erifient:
$$\sum_{i,j,[u]}\nu^{i,j}([u],[\lambda])\ell_{[u]}
(\ell_{[u]},\ell_{[\lambda]})^{-1}=
{m}_{g_s}([\lambda])\qquad (*)
$$et une \'egalit\'e de m\^eme type en \'echangeant les
r\^oles de $[u]$ et de $[\lambda]$.

 On \'ecrit le diagramme:
\[\begin{array}{ccccc} {\mathbb
C}[\hat{{\cal{W}}}_{\underline{n}}]
&\stackrel{res}{\rightarrow} &{\mathbb
C}[\hat{{\cal{W}}}_{DD_{\underline{n}}(\chi)}]
&\stackrel{ind}{\rightarrow} &{\mathbb
C}[\hat{{\cal{W}}}_{D(\chi)}]\\
\downarrow res&&\downarrow res &&\downarrow res\\ {\mathbb
C}[\hat{{\cal{W}}}_{D(\chi,g_s)}]
&\stackrel{res}{\rightarrow} &{\mathbb
C}[\hat{{\cal{W}}}_{D_{\underline{n}}(D(\chi,g_s))}]
&\stackrel{ind}{\rightarrow} &{\mathbb
C}[\hat{{\cal{W}}}_{D(\chi,g_s)}]\\ 
\downarrow ind&&\downarrow ind &&\downarrow ind\\ {\mathbb
C}[\hat{{\cal{W}}}_{D(g_s)}] &\stackrel{res}{\rightarrow}
&{\mathbb C}[\hat{{\cal{W}}}_{DD(g_s)}]
&\stackrel{ind}{\rightarrow} &{\mathbb
C}[\hat{{\cal{W}}}_{D(g_s)}]
\end{array}
\]Expliquons ce qu'est l'objet central, le
${{\cal{W}}}_{D_{\underline{n}}(D(\chi,g_s))}$ est un
sous-groupe de
${{\cal{W}}}_{D(D(\chi,g_s))}$; les collections
$\nu^{i,j}([u],[\lambda])$ qui servent \`a le construire sont
astreintes aux relations de (*) mais aussi \`a, pour tout
$i\in \{',''\}$ et pour tout $[u]\in [VP(\chi)]$:
\[\sum_{j,[\lambda]}\nu^{ij}([u],[\lambda])\ell_{[\lambda]}(
\ell_{[u]},\ell_{[\lambda]})^{-1}= n^i[u].
\] Un diagramme comme celui-ci
 est commutatif, le seul
 point est une formule \`a la Mackey du genre $res \circ ind
=ind \circ res$; une telle formule n\'ecessite des sommes:
pr\'ecis\'ement consid\'erons un groupe fini $H$ avec des
sous-groupes $H',H''$. Soit aussi une repr\'esentation de
dimension finie, $\rho'$ de $H'$ et on calcule la restriction
\`a $H''$ de l'induite de $\rho'$ \`a $H$. Cette restriction
est isomorphe \`a la somme sur $\gamma$ dans un ensemble de
repr\'esentants des doubles classes $H'\backslash H/H''$ des
induites \`a $H''$ de la repr\'esentation $\rho'$
transport\'ee par $\gamma$ et restreinte au groupe
$\gamma^{-1}H'\gamma\cap H''$ (dans \cite{mw}, ce
raisonnement est utilis\'e en 3.19 ce qui suit (4)). Comme en
loc.cit il y a la difficult\'e que les inclusions ne sont pas
compl\`etement \'evidentes. On applique cette formule 2 fois,
pour le carr\'e en bas \`a gauche, et pour le carr\'e en haut
\`a droite. Pour le carr\'e en bas \`a gauche, on l'applique
avec
$H'={\cal W}_{D(\chi,g_s)}$, $H={\cal W}_{D(g_s)}$ et
$H''={\cal W}_{DD(g_s)}$. Les doubles classes sont
pr\'ecis\'ement index\'ees par $DD(\chi,g_s)$; en effet, pour
$[\lambda]\neq [\pm 1], [\lambda] \in [VP(g_s)]$, on a \`a
consid\'erer les doubles classes:
$$
\times_{[u]\in [VP(\chi)]}\mathfrak{S}_{\nu'([u],[\lambda]}
\backslash
\mathfrak{S}_{m'([\lambda)]}/
\mathfrak{S}_{m^{','}([\lambda])}\times
\mathfrak{S}_{m^{',''}(\lambda)} $$et un objet analogue o\`u
$'$ est remplac\'e par $''$, en tenant compte du fait que
l'inclusion de $\mathfrak{S}_{\nu'([u],[\lambda])}$ dans
$\mathfrak{S}_{m'([\lambda])}$ est d\'ecrite dans ce qui
pr\'ec\`ede l'\'enonc\'e (il faut multiplier les cycles des
permutations par $\ell_{[u]}/(\ell_{[u]},\ell_{[\lambda]})$).
L'ensemble de ces doubles classes est bien index\'e par les
collections $ (\nu^{i,j}([u],[\lambda]);i,j\in {',''}, [u]\in
[VP(\chi)]$ soumises aux conditions:
$$\forall [u]\in [VP(\chi)], \forall i\in \{',''\},
\nu^{i,'}([u],[\lambda])+\nu^{i,''}([u],[\lambda])=
\nu^i([u],[\lambda]),$$
$$
\forall i,j\in \{',''\}, \sum_{[u]\in
[VP(\chi)]}\nu^{i,j}([u],[
\lambda])\ell_{[u]}/(\ell_{[u]},\ell_{[\lambda]})
=\nu^{i,j}([\lambda]).
$$
 Et ensuite il reste \`a identifier
$\times_{(\nu^{i,j}([u],[\lambda]);i,j\in {',''}, [u]\in
[VP(\chi)]}\times_{i,j \in {',''},[u]\in [VP(\chi)]}
\mathfrak{S}_{\nu^{i,j}([u],[\lambda])}$ avec
$\times_{\gamma} \gamma^{-1}H'\gamma\cap H''$ (avec les
notations pr\'ec\'edentes).  Si
$\lambda=\pm 1$, dans les objets ci-dessus, il faut remplacer
certains groupes sym\'etriques par des groupes de Weyl de
type $C$; cela ne change rien.

On a un raisonnement du m\^eme type \`a faire pour le carr\'e
en haut \`a droite du diagramme. 

Le  point maintenant \`a consid\'erer est que $\rho\circ
\iota$ n'est pas exactement
$ind\circ res$ \'ecrit sur les lignes; il faut tordre les
\'el\'ements de la colonne du milieu.  Le m\^eme
ph\'enom\`ene se produit pour
$loc_{g_s}$ et c'est ce qui motive l'introduction de
l'endomorphisme $X_{cusp}$: pour $[u]\in [VP(\chi)]$ et
$[\lambda]\in [VP(g_s)]$ regardons par quelle fonction il faut
multiplier le facteur ${\mathbb C}[\otimes_{i,j\in\{',''\}}
\hat{\cal{W}}_{\nu^{i,j}([u],[\lambda])}]$ avant d'induire
pour arriver dans ${\mathbb C}[\hat{\cal{W}}_{D(g_s)}]$ pour
que l'on obtienne le m\^eme r\'esultat qu'en faisant le
chemin, premi\`ere ligne horizontale et derni\`ere ligne
verticale. Dans toute la discussion ci-dessous, on n\'eglige,
dans les fl\`eches verticales tous les termes d\'ependant
sym\'etriquement de $\nu^{i,j}(.,.)$, sym\'etriquement en
$i,j$; c'est ce que l'on peut appeler de la torsion
sym\'etrique car elle ne g\^ene pas la commutation du
diagramme.

Supposons d'abord que $\lambda\neq \pm 1$. Si $[u]\neq \pm
1$, la ligne horizontale multiplie par le signe
$(-1)^{\nu^{'',''}([u],[\lambda])}$ et la ligne verticale
n'introduit que de la  torsion sym\'etrique; si on fait le
chemin de gauche, i.e. premi\`ere ligne verticale et
derni\`ere ligne horizontale, c'est pareil et l'on n'a pas de
probl\`eme de commutation.

Si $[u]=\pm 1$, la ligne horizontale du haut tensorise par
$sgn_{CD}w_{\nu^{'',''}([u],[\lambda])}\tilde{\chi}$  si
$\zeta_u=(-1)^{(I_u-1)/2}$; si cette \'egalit\'e n'est pas
v\'erifi\'ee c'est une autre torsion mais il faut alors aussi
tenir compte de la torsion dans la d\'efinition de $k_{\chi}$
et la combinaison des 2 ram\`enent \`a la formule
$sgn_{CD}w_{\nu^{'',''}([u],[\lambda])}\tilde{\chi}$. Quand on
fait l'autre chemin, on trouve la multiplication par
$\tilde{\chi}$ qui est introduite par l'application
$\tilde{\empty}$ puir le signe
$(-1)^{\nu^{'',''}([u],[\lambda])}$; ces 2 signes
co\"{\i}ncident gr\^ace \`a la d\'efinition de l'inclusion
donn\'ee en \ref{localisation}.

Reste le cas o\`u $\lambda=\pm 1$; on note alors 
$\epsilon'$ au lieu de $\lambda$; on rappelle que les
inclusions des groupes $\mathfrak {S}_m$ dans $W_m$ (pour $m$
un entier) consid\'er\'ees  sont telles que $(-1)^m$ est aussi
la valeur du $sgn_{CD}$ de l'image par l'inclusion. On ne
parlera donc que de $sgn_{CD}$. A priori il y a une
diff\'erence quand $[u]\neq \pm 1$ et son contraire mais comme
ci-dessus, cette diff\'erence s'efface quand on tient compte
de la d\'efintion de $k_{\chi}$; on oublie aussi le signe
$\tilde{\chi}$ qui est pris en compte par l'application
$\tilde{\empty}$. Ainsi la premi\`ere ligne horizontale et la
d\'efinition de
$k_\chi$ introduisent la multiplication par
$sgn_{CD}(w_{\nu^{'',''}([u],\epsilon')})$; la ligne
verticale multiplie par le signe de la forme
$$\prod_{i=',''}(sgn_{CD}
w_{\nu^{i,'}([u],\epsilon')})^{(I_{\epsilon_I}-1)/2}
(sgn_{CD}w_{\nu^{i,''}([u],\epsilon')})^{(I_{\epsilon_P}-1)/2}\times
\begin{cases}1\text{ si }\zeta_{\epsilon_P}=+\\
w_{\nu^{i,''}([u],\epsilon')}\text{ si
}\zeta_{\epsilon_P}=-\end{cases}
$$ et si $\epsilon'=-$ il faut encore multiplier par le
caract\`ere 
$$
\begin{cases}1 \text{ si }\epsilon_I=+
\\\prod_{i=',''}sgn_{CD} w_{\nu^{i,'}([u],-)}\text{ si
}\epsilon_I=-\end{cases}\times
\begin{cases} 1 \text{ si }\epsilon_P=+\\
\\\prod_{i=',''}sgn_{CD} w_{\nu^{i,''}([u],-)}\text{ si
}\epsilon_P=-\end{cases}\times
$$
$$
\begin{cases}1 \text{ si }(I_{\epsilon_I}+I_{\epsilon_P})/2
\text{ est impair}\\
\\prod_{i=',''}sgn_{CD} w_{\nu^{i,\delta}([u],\epsilon'}
\text{ o\`u $\delta='$ si }(I_{\epsilon_I}-1)/2 \text{ est
impair et } (I_{\epsilon_P}-1)/2 \text{ est pair et
}\delta=''\text{ sinon.}
\end{cases}
$$ Par l'autre chemin, l'application verticale introduit un
caract\`ere similaire \`a celui qui vient d'\^etre \'ecrit
sauf que ce qui
\'etait
$\nu^{i,'}$ devient
$n^{',i}$ et ce qui \'etait $\nu^{i,''}$ devient $n^{'',i}$.
Il n'y a donc pas de difficult\'e quand ce qui intervient
vraiment est un produit sur $(i,j)\in
\{',''\}$. C'est le cas quand $\epsilon_I=\epsilon_P$ et
$\zeta_{\epsilon_P}=+$. L'introduction du $X_{cusp}$ est
exactement fait pour r\'esoudre les autres cas. V\'erifions
la commutativit\'e du diagramme; on pose $\zeta=0$ si
$\zeta_{\epsilon_P}=+$ et $1$ sinon et on pose aussi
$\epsilon=0$ si $\epsilon_I=\epsilon_P$ et $1$ sinon et
finalement, on pose $\Sigma:=
(I_{\epsilon_I}+I_{\epsilon_P})/2)$. On v\'erfie que
$X_{cusp}$ n'est autre que le produit
$X_+^{1+\zeta+\Sigma}X_-^{1+\epsilon+\zeta+\Sigma}$. On
\'etudie le chemin horizontal puis vertical; il 
s'introduit donc, d'abord le signe
$sgn_{CD}w_{\nu^{'',''}([u],\epsilon')}$ puis par la
derni\`ere fl\`eche verticale, un signe $\chi_{cusp}$. Mais
pour les probl\`emes de commutation, on peut multiplier ce
signe par n'importe quel  signe de la forme
$\prod_{i,j\in
\{',''\}}sgn_{CD}w_{\nu^{i,j}([u],\epsilon'_0)}$, o\`u
$\epsilon'_0\in \{\pm 1\}$ comme expliqu\'e ci-dessus. Ce qui
veut dire qu'au lieu d'utiliser $\chi_{cusp}$ tel qu'il a
\'et\'e
\'ecrit, on peut utiliser 
$$
(\chi''_+\chi''_-)^{1+\Sigma+\delta}(\chi''_-)^\epsilon=
(\chi''_+)^{1+\delta+\Sigma}(\chi''_-)^{1+\delta+\epsilon+
\Sigma}
.\eqno
(*)
$$
Quand on fait la derni\`ere fl\`eche verticale, en terme de
$w_{\nu^{i,j}}$ cela devient un produit sur tout $[u]\in
[VP(\chi)]$ de
$$
\prod_{i\in\{',''\}}sgn_{CD}w_{\nu^{i,''}([u],+)}
^{1+\delta+\Sigma}\prod_{i\in\{',''\}}
sgn_{CD}w_{\nu^{i,''}([u],-)}
^{1+\delta+\epsilon+\Sigma}.$$
En incorporant le signe de la ligne horizontale, on trouve,
un produit sur tout $[u]$
$$
sgn_{CD}w_{\nu^{',''}([u],+)}
^{1+\delta+\Sigma}sgn_{CD}w_{\nu^{'',''}([u],+)}
^{\delta+\Sigma}sgn_{CD}w_{\nu^{',''}([u],-)}
^{1+\delta+\epsilon+\Sigma}
sgn_{CD}w_{\nu^{'',''}([u],-)}
^{\delta+\epsilon+\Sigma}. \eqno (**)
$$
On examine maitenant le chemin utilisant d'abord la
premi\`ere fl\`eche verticale puis l'action de $X_{cusp}$ et
la derni\`ere ligne horizontale. Dans $X_{cusp}$ on commence
par multiplier par un caract\`ere qui est exactement le
caract\`ere (*) qui s'introduit par la fl\`eche verticale 
 apr\`es la simplification  effectu\'ee ci-dessus. Finalement,
pour ce chemin, il suffit de regarder le caract\`ere de la
derni\`ere ligne horizontale en tenant compte de l'inversion
\'eventuelle. Or on a inversion entre $m'(+)$ et $m''(+)$ par
hypoth\`ese si $\delta+\Sigma$ est impaire et inversion entre
$m'(-)$ et $m''(-)$ si $\delta+\epsilon+\Sigma$ est impaire.
L'inversion
entre $'$ et $''$ a pour effet que $\rho\circ \iota$ introduit
le signe
$sgn_{CD}w_{\nu^{',''}([u],\epsilon')}$ au lieu de $sgn_{CD}
w_{\nu^{'',''}([u],\epsilon')}$.
On trouve donc exactement le caract\`ere (**). Cela termine
la preuve.

\section{Stabilit\'e}

\subsection {Stabilit\'e, d\'efinition\label{stabilite}}

On reprend encore les d\'efinitions de \cite{mw}; soit $G$
un groupe classique qui est donc le groupe des automorphismes
d'une forme (ici orthogonale ou unitaire); on doit
consid\'erer simultan\'ement 2 formes de ce groupe
correspondant \`a 2 formes orthogonale ou unitaire,
s\'epar\'ee dans le cas orthogonal par l'invariant de Hasse
et dans le cas unitaire par la parit\'e de la valuation du
d\'eterminant. On note ces 2 formes
$G_{iso}$ et
$G_{an}$ en imposant que $G_{iso}$ est la forme
quasid\'eploy\'ee et que $G_{an}$ est l'autre groupe; dans
tous les cas, $G_{an}$ est une forme int\'erieure de
$G_{iso}$ mais \'eventuellement, on a m\^eme un isomorphisme
$G_{an}\simeq G_{iso}$; ces 2 formes interviennent dans le
calcul du centralisateur d'un \'el\'ement semi-simple tel que
fait dans \ref{localisation} et les constructions d\'ependent
de la forme orthogonale ou unitaire qui intervient et pas
seulement de son groupe d'automorphismes, d'o\`u la
n\'ecessit\'e de garder la diff\'erence dans les notations. On
sait d\'efinir la classe de conjugaison stable de tout
\'el\'ement fortement r\'egulier de $G_{\sharp}$ pour
$\sharp=iso $ ou $an$ et on sait aussi d\'efinir une
inclusion de l'ensemble des classes de conjugaison stable de
$G_{an}$ dans  l'ensemble des classes de conjugaison stable
dans $G_{iso}$. Soit $\phi=(\phi_{iso},\phi_{an})$ une
fonction dans
$C^{\infty}_c(G_{iso})\oplus C^{\infty}_c(G_{an})$; on dit
qu'elle est stable si les int\'egrales orbitales de
$\phi_{iso}$ et de $\phi_{an}$ sont constantes sur les
classes de conjugaison stable et si les int\'egrales
orbitales de
$\phi_{iso}$ et de $\phi_{an}$  se correspondent pour
l'inclusion des classes stables pour $G_{an}$ dans les
classes stables de $G_{iso}$ et 
$\phi_{iso}$ a une int\'egrale nulle sur les classes stables
de $G_{iso}$ ne provenant pas de $G_{an}$); il a
\'evidemment fallu fixer des mesures coh\'erentes. On dit que
$\phi$ est semi-stable si $\phi_{iso}$ et $\phi_{an}$ sont
stables mais si  pour tout
$\gamma$ fortement r\'egulier dans
$G_{an}$, l'int\'egrale orbitale de $\phi_{an}$ sur la classe
de conjugaison stable de $\gamma$ est l'oppos\'ee de
l'int\'egrale de $\phi_{iso}$ sur  la classe de conjugaison
stable  dans $G_{iso}$ correspondant
\`a celle de
$\gamma$. On dit que
$\phi_{iso}$ est instable si pour tout $\gamma$ fortement
r\'egulier  l'int\'egrale sur la classe de conjugaison stable
de $\gamma$ est nulle; on d\'efinit de m\^eme
$\phi_{an}$ instable et on dit que $\phi$ est instable si
$\phi_{iso}$ et $\phi_{an}$ sont instables. 

\noindent
Soit $\underline{n}\in D(\chi)$ avec une donn\'ee cuspidale
${cusp}$; on dit que
\begin{enumerate}
\item
$\underline{n}, cusp$ est stable si  
$\epsilon_I=\epsilon_P$, $\zeta_+=\zeta_-=+$,
$\vert I_\epsilon-P_\epsilon\vert=1$ pour $\epsilon=\pm$ et
$n''[u]=0$ pour tout
$[u]\in [VP(\chi)]$.
\item
On dit que $\underline{n},{cusp}$ est semi-stable si
$\epsilon_I=\epsilon_P$, $\zeta_+=\zeta_-=-$,
$\vert I_\epsilon-P_\epsilon\vert=1$ pour $\epsilon=\pm$ et
$n'([u])=0$ pour tout $[u]\in [VP(\chi)]$.
\item
On dit que $\underline{n}, {cusp}$ est instable dans tous les
autres cas.
\end{enumerate}

\noindent 
\bf Remarque: \sl soit $(\psi,\epsilon)$ un param\`etre
discret de niveau z\'ero. Le couple
$\underline{n}_{\psi,\epsilon}, cusp$ qui lui est associ\'e
avec la repr\'esentation de Springer-Lusztig est stable si et
seulement si  pour tout
$[u]\in [VP(\chi)] \neq [\pm 1]$, $\epsilon_{[u]}$ est le
caract\`ere trivial et si $U_{[\pm 1],-}=\emptyset$ (avec les
notations de \ref{classification});
$\underline{n}_{\psi,\epsilon}, cusp$  est semi-stable si
$\epsilon_{[u]}\equiv -1$ pour tout $[u]\in [VP(\chi)]\neq
[\pm 1]$ et si
$U_{[\pm 1],+}=\emptyset$. Et
$\underline{n}_{\psi,\epsilon},cusp$ est instable dans tous
les autres cas.\rm

\

\noindent
On rappelle les formules donn\'ees dans
\ref{springerlusztig}. Pour $[u]\neq [\pm 1]$,   la
traduction de $n'([u])=0$ ou $n''([u])=0$ en terme du
caract\`ere du groupe des composantes est claire .

Le cas de $u=\pm 1$  est plus
compliqu\'e. On regarde d'abord la partie cuspidale; \`a
chaque orbite $U_{u,\epsilon'}$ munie de son
caract\`ere du groupe des composantes est associ\'e un entier
$k_{u,\epsilon'}$ par la repr\'esentation de Springer
g\'en\'eralis\'ee. On fixe $u,\epsilon'\in \{\pm 1\}$
et on montre d'abord l'\'equivalence:
$$ k_{u,\epsilon'}=0\Leftrightarrow \vert
I_u-P_u\vert=1 \text{ et }
\zeta_{u}=\epsilon'.
$$ En effet, on v\'erifie d'apr\`es les formules donn\'ees que
$k_{u,\zeta_u}=(I_u+P_u-1)/2$ et
$k_{u,-\zeta_u}=(\vert
I_u-P_u\vert-1)/2$. Et l'\'equivalence est
alors claire, en tenant compte du fait que $I_u$ est
impair alors que $P_u$ est pair par hypoth\`ese.
Ensuite, c'est presque les d\'efinitions que
$U_{u,\epsilon'}=0$ est \'equivalente \`a
$k_{u,\epsilon'}=0$ et
$n^\delta(u)=0$ o\`u
$\delta='$ si $\epsilon'=+$ et $''$ si $\epsilon'=-$. 

\subsection{Stabilit\'e, th\'eor\`eme\label{theoreme}}

On fixe une donn\'ee cuspidale $cusp$ et $\underline{n}\in
D(\chi)$.

\bf Th\'eor\`eme: \sl  soit $\phi \in {\mathbb
C}[\hat{{\cal{W}}}_{\underline{n},cusp}]$ et soit
$\Phi:=k_\chi\,
\rho\circ \iota (\phi)$. Alors $\Phi$ est stable si et
seulement si
$\underline{n}, cusp$  est stable; de m\^eme  $\Phi$ est
semi-stable si et
seulement si 
$\underline{n}, cusp$ est semi-stable
et 
$\Phi$ est instable si et
seulement si 
$\underline{n}, cusp$ est instable.
\rm

On suit la m\'ethode de \cite{mw}  3.20 (qui d\'emontre le
m\^eme th\'eor\`eme dans le cas o\`u
$[VP(\chi)]=[1]$. On \'ecrit $\Phi:=(\Phi_{iso},\Phi_{an})$.
On fixe un \'el\'ement semi-simple fortement r\'egulier
$g\in SO(2n+1,F)_{iso}$ et on \'etudie les int\'egrales
orbitales de $\Phi_{iso}$ pour les \'el\'ements de la classe
de conjugaison stable de $g$ ainsi que celles de
$\Phi_{an}$ pour la classe de conjugaison stable dans
$SO(2n+1,F)_{an}$ quand elle existe. Il est clair que ces
int\'egrales orbitales sont nulles si $g$ n'est pas
elliptique et compact. On \'ecrit $g=g_sg_u$ comme en
\ref{localisation}. L'ensemble
$[VP(g_s)]$ est ind\'ependant de $g$ dans sa classe de
conjugaison stable et quand $g$ varie dans sa classe de
conjugaison stable vue dans $SO(2n+1,F)_{iso}\cup
SO(2n+1,F)_{an}$
$g_s$ varie exactement dans sa classe de conjugaison stable
dans $SO(2n+1,F)_{iso}\cup SO(2n+1,F)_{an}$. Les classes de
conjugaison dans la classe de conjugaison stable de $g_s$
sont param\'etr\'ees (\cite{w1} 1.7) par les collections
$\{\sharp_{[\lambda]}\in \{+1,-1\}\simeq
\{iso,an\}\}_{[\lambda]\in [VP(g_s)]}$ de telle sorte que si
$g_s$ correspond \`a la collection 
$\{\sharp_{[\lambda]}(g_s); \lambda\in [VP(g_s)]\}$
 le commutant de $g_s$ est isomorphe au produit
$Aut((F'_{[\lambda]})^{m_{[\lambda]}},<\, ,\,
>_{\sharp_{[\lambda]}(g_s)})$ o\`u 
$F'_{[\lambda]}$ est une extension non
ramifi\'ee de degr\'e 2 de $F_{[\lambda]}$ l'extension
non ramifi\'ee de
$F$ de degr\'e
$\ell_{[\lambda]}$ ($2\ell_{[\lambda]}$ est le
cardinal de l'ensemble $[\lambda]$) et o\`u $<\, ,\,
>_{\sharp_{[\lambda]}(g_s)}$ est une forme unitaire (pour
l'extension $F'_{[\lambda]}$ de $F_{[\lambda]}$) dont le
d\'eterminant est de valuation paire ou impaire suivant que
$\sharp_{[\lambda]}(g_s)=1$ ou $-1$, si $[\lambda]\neq \pm 1$
et est une forme orthogonale si $[\lambda]=[\pm 1]$ (il n'y a
alors pas d'extension de degr\'e 2 \`a consid\'erer); dans ce
dernier cas, on a la m\^eme propri\'et\'e que
pr\'ec\'edemment mais ''parit\'e de la valuation du
d\'eterminant'' \'etant remplac\'e par invariant de Hasse.
Pour d\'ecrire les classes de conjugaison dans la classe de
conjugaison stable de $g$, il faut encore d\'ecrire o\`u
varie $g_u$ quand $g_s$ est fix\'e. Comme on appliquera le
d\'ebut de la section 3 de
\cite{mw} tel quel, nous n'avons pas besoin de faire cette
description et on renvoie le lecteur \`a loc. cite.

Fixons maintenant $g=g_sg_u\in SO(2n+1,F)_{iso}\cup
SO(2n+1,F)_{an}$ comme ci-dessus et calculons
$I_{g}(\Phi)$. On note $\sharp(g)$ l'\'el\'ement $iso$ ou
$an$ tel que $g\in SO(2n+1,F)_{\sharp(g)}$.  On a d\'efini la
fonction de Green
$Q(loc_{g_s}(\Phi))$ en \ref{localisation}; c'est une fonction
\`a support les \'el\'ements topologiquement unipotents  sur
le groupe
$$\times_{[\lambda]\in [VP(g_s)], [\lambda]\neq
\pm
1}U(m(\lambda),F'_{[\lambda]}/F_{[\lambda]})_{\sharp_{[\lambda]}(g_s)}\times
SO(m([1]),F)_{\sharp_{[1]}(g_s)}\times
O(m([-1]),F)_{\sharp_{[-1]}(g_s)},$$ 
o\`u les notations sont
celles de \ref{localisation}.  L'\'el\'ement
$g_u$ d\'efinit une classe de conjugaison d'\'el\'ements
topologiquement unipotents dans ce groupe. Pour la suite on
notera
$(loc_{g_s}(\Phi))_{[\lambda]}$ la fonction sur le groupe
index\'e par $[\lambda]$ d\'efinie par
$(loc_{g_s}(\Phi))$ quand les points dans les  groupes
index\'es par $[\lambda']\neq [\lambda]$ sont fix\'es.

Fixons une classe de conjugaison stable dans $SO(2n+1,F)$
d'\'el\'ements semi-simples r\'eguliers ${\cal C}_{st}$; sans
restreindre la g\'en\'eralit\'e, on les suppose compacts
(sinon les int\'egrales orbitales sont nulles). On note
g\'en\'eriquement ${\cal C}$ les classes de conjugaison
incluses dans ${\cal C}_{st}$; les classes de conjugaison le
sont pour un groupe, c'est-\`a-dire qu'une telle classe
${\cal C}$ correspond \`a une valeur de $\sharp$ qui est
not\'ee $\sharp({\cal C})$. Pour chaque classe ${\cal C}$, on
fixe un \'el\'ement $g({\cal C})\in {\cal C}$. On a \`a
calculer pour $\Phi$ comme ci-dessus et pour $\sharp$ fix\'e
$I^{st,\sharp}({\cal C}_{st},\Phi):=\sum_{{\cal C}\in {\cal
C}_{st},
\sharp({\cal C})=\sharp} I^{SO(2n+1,F)_{\sharp}}(g({\cal C}),
\Phi)$. On \'ecrit chaque $g({\cal C})=g_s({\cal C})g_u({\cal
C})$. On \'etablit une relation d'\'equivalence entre les 
${\cal C}\in {\cal C}_{st}$ par ${\cal C}\sim {\cal C}'$ si
$g_s({\cal C})$ est conjugu\'e de $g_s({\cal C}')$; on les
supposera alors \'egaux.  On
\'ecrira donc
$g_s([{\cal C}])$ plut\^ot que $g_s({\cal C})$.

Il existe une classe stable d'\'el\'ements semi-simples dont
les valeurs propres sont des racines de l'unit\'e d'ordre
premier
\`a
$p$, ${\cal C}_{s,st}$ tel que
$g_s([{\cal C}])$ soit un repr\'esentant des classes de
conjugaison dans ${\cal C}_{s,st}$. On l'utilisera
plus bas mais tout d'abord, on \'ecrit:
$$ I^{st,\sharp}({\cal C}_{st},\Phi):=\sum_{[{\cal C}]\in
{\cal C}_{st}/\sim,
\sharp({\cal C})=\sharp}\sum_{{\cal C}\in [{\cal
C}]}I^{Cent_{SO(2n+1,F)_{\sharp}(g_s([{\cal C}]))}}(g_u({\cal
C}),loc_{g_s([{\cal C}])}\rho\circ
\iota\, \phi).$$ On utilise \ref{cle} pour r\'ecrire, pour
$[{\cal C}]$ fix\'ee:
$$
\sum_{{\cal C}\in [{\cal
C}]}I^{Cent_{SO(2n+1,F)_{\sharp}(g_s([{\cal C}]))}}(g_u({\cal
C}),loc_{g_s([{\cal C}])}\rho\circ
\iota\, \phi)=
\sum_{{\cal C}\in [{\cal
C}]}I^{Cent_{SO(2n+1,F)_{\sharp}(g_s([{\cal C}]))}}(g_u({\cal
C}),Q\,\rho\circ
\iota\,X_{cusp} loc_{g_s([{\cal C}])}\tilde{\empty}
(\phi)).$$On utilise  tout de suite le fait que
$[VP(g_s([{\cal C}]))]$ est ind\'ependant de
$[{\cal C}]$ dans ${\cal C}_{st}$; ce qui varie sont  les
invariants  des formes $<\, ,\,>_{[\lambda]}$, cf.
\ref{localisation} et ci-dessus. On  remplace donc la notation
$[VP(g_s)]$ par $[VP({\cal C}_{st})]$. On d\'ecompose la
somme ci-dessus en produit sur les
$[\lambda]\in [VP({\cal C}_{st})]$ et on constate qu'elle est
nulle si l'une des  composantes
$Q\, \rho\circ \iota X_{cusp}(loc_{g_s([{\cal
C}]}\tilde{\empty}(\phi))_{[\lambda]}$ est instable.  La
condition d'instabilit\'e pour ce genre de fonction (c'est
\`a dire pour des fonctions dans l'image de $Q\,
\rho\circ\iota$) est d\'ecrite dans
\cite{mw} 3.4 pour les groupes unitaires (c'est-\`a-dire ici
pour  $[\lambda]\neq \pm 1$) et en
\cite{mw}  3.12 pour les groupes orthogonaux (c'est-\`a-dire
ici pour 
$[\lambda]=[\pm 1]$). Pour
$[\lambda]\neq [\pm 1]$, on a instabilit\'e si
$m'([\lambda])m''([\lambda])\neq 0$. Pour $\lambda\in \{\pm
1\}$, on a  les donn\'ees pour le support cuspidal des
fonctions de Green qui sont
$\vert
\tilde{cusp}\vert$ c'est-\`a-dire:
$\vert r'_\epsilon\vert:=(I_++\epsilon\delta I_-)/2$, o\`u
$\delta=+$ ou
$-$ de fa\c con \`a ce que $\vert r '_+\vert$ soit impair et
$\vert r ''_\epsilon\vert:=\vert \bigr(
(-1)^{(I_+-1)/2}\zeta_{+}P_++\epsilon(-1)^{(I_--1)/2}\zeta_{-}P_-
\bigr)/2\vert$. Ainsi $\delta=(-1)^{1+(I_++I_-)/2}$ et $\vert
r''_\epsilon\vert=\vert
(P_++\epsilon\delta\zeta_+\zeta_-P_-)\vert.
$

 Pour que les int\'egrales ne soient pas
nulles, il faut que $M'(\lambda):=m'(\lambda)-
(r'_{\lambda})^2$ soit un entier pair ($\geq 0$) et la m\^eme
propri\'et\'e pour 
$M''(\lambda):=m''(\lambda)- (r''_{\lambda})^2$. Pour
$\epsilon=\pm$, le terme correspondant \`a $\lambda=\epsilon$ 
est instable si et seulement si soit $\vert\vert
r'_\epsilon\vert-\vert r''_\epsilon\vert\vert > 1$ soit
$r'_\epsilon r''_\epsilon M''(\epsilon)$ ou soit
$M'(\epsilon)M''(\epsilon)=0$. On remarque que

$$\vert r'_\epsilon\vert-\vert r''_\epsilon\vert=
\vert I_++\epsilon\delta I_-\vert-\vert P_++\epsilon\delta
\zeta_+\zeta_-P_-\vert.\eqno(*)$$ 
 
On a donc instabilit\'e si l'une des conditions $ I_++I_-
-\vert P_++\zeta_+\zeta_- P_-\vert \in \{-2,0,2\}$, $ \vert
I_--I_-\vert  -\vert P_+-\zeta_+\zeta_- P_-\vert \in
\{-2,0,2\}$ n'est pas satisfaite. On remarque
que $I_++I_--\vert I_+-I_-\vert \geq 2$ et que si $P_+P_-\neq
0$, $P_++P_--\vert P_+-P_-\vert\geq 4$ pour des questions de
parit\'e. Ainsi si
$\zeta_+\zeta_- \neq +$, c'est-\`a-dire vaut - et si 
$P_+P_-\neq 0$ la diff\'erence entre
$I_++I_--\vert P_+-P_-\vert$ et $\vert I_+-I_-\vert -(
P_++P_-)$ est au moins 6. On ne peut donc avoir les deux
conditions satisfaites  en m\^eme temps. Ainsi, on a
instabilit\'e si 
$P_+P_-\neq 0$ mais
$\zeta_+\zeta_-=-$. Si $P_+P_-=0$, en reprenant les
notations, $\epsilon_I$ et
$\epsilon_P$ de \ref{localisation}, on a donc
$P_{-\epsilon_P}=0$ et on v\'erifie que n\'ecessairement
$I_{-\epsilon_I}=1$.

On a donc d\'ej\`a d\'emontr\'e que l'on a instabilit\'e sauf
si soit
$\zeta_+\zeta_-=+$ soit $P_{-\epsilon_P}=0$ et
$I_{-\epsilon_I}=1$.

R\'ecrivons les conditions (*) ci-dessus sous la forme plus
simple
$(I_{\epsilon_I}-P_{\epsilon_P})\pm
(I_{-\epsilon_I}-P_{-\epsilon_P})\in \{-2,0,2\}$ et encore
$$(I_{\epsilon_I}-P_{\epsilon_P})\in
\{-1,1\}; (I_{-\epsilon_I}-P_{-\epsilon_P})\in
\{-1,1\}.\eqno(*)_{cusp}$$
On rappelle la convention que $\epsilon_I=\epsilon_P$ si
$(I_+-I_-)(P_+-P_-)=0$ et que $\epsilon_I=\epsilon_P=+$ si
$I_+=I_-$ et $P_+=P_-$.

 On a donc instabilit\'e au moins
s'il existe $[\lambda]\in [VP({\cal C}_{st})]$ tel que
$m'([\lambda])m''([\lambda])\neq 0$ (en rempla\c cant $m'$
par $M'$ et $m''$ par $M''$ si
$\lambda=\pm 1$ ou si (*)$_{cusp}$ n'est pas satisfaite ou
encore s'il existe
$\epsilon=\pm $ tel que
$M'([\epsilon])=0$,
$M''(\epsilon)\neq 0$ et $r'_\epsilon r''_\epsilon \neq 0$.

Par les r\'ef\'erences d\'ej\`a donn\'ees, on sait aussi
quand ces sommes partielles ne d\'ependent que de l'invariant
$\sharp_{[\lambda]}$ de $<\, ,\, >_{[\lambda]}$ ou sont
ind\'ependantes du choix de $[{\cal C}]$ dans ${\cal
C}_{g_s}$. Il faut simplement faire attention que
$loc_{g_s{[{\cal C}]}}$ d\'epend de $[{\cal C}]$ dans ${\cal
C}_{st}$ par ce qui est not\'e $c_{cusp}(g_s)$ en
\ref{localisation} et donc (en supprimant ce qui est encore
ind\'ependant) par le signe:
$$(\eta'_+\eta'_-)^{(I_{\epsilon_I}-1)/2}(\eta''_+\eta''_-)^{
(I_{\epsilon_P}-1)/2}\times \begin{cases} 1\text{ si
}\zeta_{\epsilon_P}=+\\
\eta''_+\eta''_- \text{ si } \zeta_{\epsilon_P}=- \end{cases}
\times \begin{cases} 1\text{ si }\epsilon_I=\epsilon_P=+ \\
\eta'_-\eta''_-\text{ si } \epsilon_I=\epsilon_P=-\\
\eta'_-\text{ si } \epsilon_I=-,\epsilon_P=+\\
\eta''_- \text{ si } \epsilon_I=+,\epsilon_P=-.\end{cases}
$$ Il est temps d'utiliser les propri\'et\'es des formes $<\,
,\, >_{\epsilon}$ pour $\epsilon=\pm$ et de leur
``r\'eduction'' (rappel\'ees en \cite{mw} 3.11):
$\eta'_\epsilon\eta''_\epsilon$ est calcul\'e par le
discriminant de la forme $<\, ,\, >_\epsilon$ et
$\eta'_\epsilon$ comme $\eta''_\epsilon$ sont soit
l'invariant de Hasse soit son oppos\'e (le choix d\'epend du
discriminant). Ainsi, dans la formule ci-dessus, si
$\epsilon_I\neq \epsilon_P$, $c_{cusp}(g_s([{\cal C}]))$
d\'epend de l'invariant de Hasse de l'une des formes $<\, ,\,
>_\epsilon$ et pas de l'autre: en effet le premier terme
d\'epend du produit des 2 invariants de Hasse, le deuxi\`eme
terme d\'epend soit du produit soit vaut 1 et le troisi\`eme
terme d\'epend de l'un des invariants de Hasse exactement. Par
contre si $\epsilon_I=\epsilon_P$ alors $c_{cusp}({\cal C})$
d\'epend du produit des invariants de Hasse quand
$\zeta_{\epsilon_P}=-$ et est constant si
$\zeta_{\epsilon_P}=+$.

 Comme la seule chose qui est fix\'ee pour les classes de
conjugaison dans la classe de conjugaison stable de $g_s$
\`a l'int\'erieur d'un groupe
$SO(2n+1,F)_{\sharp}$ o\`u $\sharp$ est fix\'e est le produit
sur tous les $[\lambda]$ des invariants $\sharp_{[\lambda]}$
(qui sont les invariants de Hasse pour $\lambda=\pm 1$),  on
voit ais\'ement que l'on a encore instabilit\'e s'il existe
$[\lambda]\in [VP({\cal C}_{st})]$ tel que le terme
correspondant d\'epend de l'invariant $\sharp_{[\lambda]}$ de
la forme
$<\, ,\, >_{[\lambda]}$ et qu'il existe $\lambda'\in 
[VP({\cal C}_{st}]$ tel que le terme correspondant ne
d\'epend pas de l'invariant  de la forme
$<\, ,\, >_{[\lambda']}$. Et on aura stabilit\'e si aucun des
termes n'en d\'epend et semi-stabilit\'e si tous les termes
en d\'ependent. Ainsi la stabilit\'e se produit quand
$m''([\lambda])=0$ pour tout $[\lambda] \in [VP({\cal
C}_{st})]-\{-1,+1\}$ et si le produit des termes correspondant
\`a $+1$ et $-1$ est aussi ind\'ependant des invariants de
Hasse. Comme on l'a vu ci-dessus, il faut distinguer le cas
$\epsilon_I\neq
\epsilon_P$ du cas o\`u l'on a \'egalit\'e. Supposons d'abord
que $\epsilon_I=\epsilon_P$; dans ce cas, si
$\zeta_{\epsilon_P}=+$, $c_{cusp}({\cal C})$ est ind\'ependant
des invariants de Hasse, il faut donc aussi que les
int\'egrales en soient ind\'ependantes et donc que 
$M''(+1)=M''(-1)=0$.

Par contre si $\zeta_{\epsilon_P}=-$, toujours sous
l'hypoth\`ese $\epsilon_I=\epsilon_P$, il faut $\vert
r'_\epsilon\vert 
\vert r''_\epsilon\vert =0$ et
$M'(\epsilon)=0$ pour
$\epsilon=\pm$. Pour avoir stabilit\'e, on a d\'ej\`a
vu qu'il faut $P_{-\epsilon_P}=0$ et $I_{-\epsilon_I}=1$. La
condition
$\vert
r'_\epsilon\vert 
\vert r''_\epsilon\vert =0$
\'ecrite pour $\epsilon=\delta$, donne 
$(I_++I_-)\vert (P_+--P_-)\vert=0$. D'o\`u $P_+=P_-=0$ et on
retrouve aussi $I_{\epsilon_I}=1$ en utilisant $(*)_{cusp}$.
D'o\`u par convention $\zeta_+=\zeta_-=+$ ce qui contredit
$\zeta_{\epsilon_P}=-$.

Terminons le cas de la stabilit\'e quand
$\epsilon_I=\epsilon_P$; la stabilit\'e est alors
\'equivalente \`a $m''([\lambda])=0$ pour tout $[\lambda]\neq
\pm 1$ et pour $\lambda=\pm 1$, il faut $M''(\lambda)=0$,
$\zeta_+=\zeta_-=+$ et (*)$_{cusp}$. On remarque que la
condition
$\epsilon_I=\epsilon_P$ coupl\'ee avec (*)$_{cusp}$ est
\'equivalente \`a $\vert I_+-P_+\vert=1$ et $\vert
I_--P_-\vert=1$. Il faut encore utiliser le fait que la
localisation est non nulle si l'ensemble $D(\chi,g_s)$ n'est
pas vide; c'est-\`a-dire qu'il existe donc une collection
$\nu'([u],[\lambda]),\nu''([u],[\lambda])$ satisfaisant \`a:
$$ m'([\lambda])=\sum_{[u]\in [VP(\chi)]}
\nu'([u],[\lambda])\ell_{[u]}/(\ell_{[u]},(\ell_{[\lambda]}),
\qquad m''([\lambda] =
\sum_{[u]\in [VP(\chi)]}
\nu''([u],[\lambda])\ell_{[u]}/(\ell_{[u]},(\ell_{[\lambda]})
,$$
pour tout $[\lambda]\neq [\pm 1]$ et une formule analogue
quand $\lambda\in \{\pm 1\}$. On a aussi, avec les m\^emes
notations, pour tout
$[u]\in [VP(\chi)]\neq \pm 1$:
$$ n'([u])=\sum_{[\lambda]\in [VP([{\cal
C}])]}\nu'([u],[\lambda])\ell_{[\lambda]}
/(\ell_{[u]},(\ell_{[\lambda]})
,\qquad
\nu''([u])=\sum_{[\lambda]
\in [VP([{\cal C}])]}\nu''([u],[\lambda]) \ell_{[\lambda]}
/(\ell_{[u]},(\ell_{[\lambda]});
$$et une formule analogue pour $u\in \{\pm 1\}$. On en
d\'eduit que les conditions
$n''([u])=0$ (si
$u\neq \pm 1$) et
$N''(\pm 1)=0$ sont \'equivalentes \`a leurs analogues pour
$m',m''$ et $[\lambda]$.

Pour la stabilit\'e, reste \`a voir le cas o\`u
$\epsilon_I\neq
\epsilon_P$. On a vu que $c_{cusp}(g_s({\cal C}))$ d\'epend de
l'un des invariants de Hasse; il faut donc que les
int\'egrales d\'ependent elles aussi de l'un des invariants
de Hasse exactement. Mais on doit donc avoir l'une des
conditions
$r'_\epsilon r''_\epsilon=0$ qui n\'ecessairement
entra\^{\i}ne soit $I_+=I_-$ soit $P_+=P_-$. Et on a donc
imm\'ediatement une impossibilit\'e avec les conventions sur
$\epsilon_I$ et $\epsilon_P$.

Cela termine la preuve en ce qui concerne la stabilit\'e.

Pour la semi-stabilit\'e: le raisonnement est du m\^eme type,
il faut pour tout $\lambda\neq \pm 1$, $m'(\lambda)=0$. Pour
$\lambda=\pm 1$, on v\'erifie qu'il faut
$\epsilon_I=\epsilon_P$ (c'est comme ci-dessus), puis
$M''(\lambda)=0$ et $\zeta_+=\zeta_-=-$. Ensuite, on se
rappelle des \'echanges  induits par $X_{cusp}$; on doit
\'echanger $\nu'([u],\epsilon')$ et $\nu''([u],\epsilon')$
pour tout $[u]\in [VP(\chi)]$ et tout $\epsilon'=\pm 1$ (cf.
\ref{localisation}). On en d\'eduit que la semi-stabilit\'e
est \'equivalente \`a ce que $n'([u])=0$ pour tout $[u]\neq
\pm 1$ dans $[VP(\chi)]$,
$n'(\pm 1)=0$ ainsi que les conditions d\'ej\`a \'ecrites sur
la partie cuspidale. Cela termine la preuve.

\subsection{Traduction en termes de
param\`etres\label{traduction}}

La remarque ci-dessous vient d'id\'ees de Lusztig avec des
compl\'ements pour les groupes non connexes de \cite{w3}.
On consid\`ere l'ensemble des quadruplets d'entiers positifs
ou nuls introduit dans \ref{springerlusztig}, pour
$u=\pm$ et
$\epsilon'=\pm$,
$k_{u,\epsilon'}$ et on pose $I_\epsilon$ l'entier
impair du couple $(k_{u,+}+k_{u,-}+1,\vert
k_{u,+}-k_{u,-}\vert)$ et $P_\epsilon$ l'entier
pair. On pose aussi $\zeta_u$ le signe de
$k_{u,+}-k_{u,-}$ avec la convention que si ce
nombre est nul $\zeta_u=(-1)^{k_{u,+}}$ ce qui
est compatible avec la convention de \ref{localisation} car
$(I_{u}-1)/2=k_{u,+}$ dans ce cas.

\bf Remarque: \sl avec les notations pr\'ec\'edentes, on a
l'\'equivalence des conditions:
$$
\forall u\in \{\pm\}; \vert
I_u-P_u\vert=1 \text{ et } \zeta_u=+
\Leftrightarrow \forall u\in \{\pm \}\, k_{u,-}=0;
$$
$$
\forall u\in \{\pm\}; \vert
I_u-P_u\vert=1 \text{ et } \zeta_u=-
\Leftrightarrow \forall u\in \{\pm \}\, k_{u,+}=0.
$$\rm On a pour $u=\pm$, $\vert
I_u-P_u\vert=1+k_{u,-\zeta_\epsilon}$ et la
remarque en d\'ecoule.

La remarque pr\'ec\'edente motive la d\'efinition:

\bf D\'efinition: \sl soit $\psi,\epsilon$ un param\`etre
discret de niveau 0. On reprend les notations de
\ref{classification} et \ref{springerlusztig}. On dit qu'il
est stable si pour tout
$[u]
\in [VP(\chi)]$  l'orbite
$U'_{[u]}=0$; il est semi-stable si pour tout $[u]\in
[VP(\chi)$  l'orbite $U''_{[u]}=0$ et il est instable
sinon.\rm

\bf Corollaire: \sl l'espace des fonctions associ\'ees via la
repr\'esentation de Springer-Lusztig et les faisceaux
caract\`eres \`a un param\`etre discret de niveau z\'ero est
form\'e de fonctions stables, semi-stables ou instables si et
seulement si le param\`etre est stable, semi-stable ou
instable.\rm

Avec la remarque, cela r\'esulte de \ref{theoreme} et de la
d\'efinition ci-dessus.

\section{Interpr\'etation}

\subsection{Transformation de Fourier\label{fourier}}

Dans les conjectures de Langlands, les propri\'et\'es de
stabilit\'e ne s'expriment pas comme dans le corollaire
ci-dessus. Ce ne sont pas certains param\`etres qui sont
stables  mais au contraire ce sont des combinaisons
lin\'eaires. Les 2 fa\c cons d'exprimer le r\'esultat  se
d\'eduisent l'une de l'autre par une transformation style
transformation de Fourier. Plus  pr\'ecis\'ement, on fixe
$\chi$ et on fixe des orbites $U_{[u]}$ pour tout $[u]\in
[VP(\chi)]$ v\'erifiant les conditions de \ref{classification}
et on note
$\underline{U}:=\{U_{[u]};[u]\in [VP(\chi)]\}$. On note
$\mathfrak{P}_{\chi,\underline{U}}$ l'espace vectoriel
complexe de base l'ensemble des param\`etres discrets de
niveau 0 tels que la restriction de $\psi$ \`a $I_F\times
SL(2,{\mathbb C})$ (cf. \ref{classification}) soit
d\'etermin\'ee par
$\chi$ et $\underline{U}$. Le principe est de d\'efinir un
produit scalaire sur cet espace,
$<\, ,\, >$ et de d\'efinir la transformation $\cal F$, en
posant:
\[\forall p\in \mathfrak{P}_{\chi,\underline{U}}; {\cal{F}}
(p):= \sum_{p'\in \mathfrak{P}_{\chi,\underline{U}}} <p,p'>p'
\]Et si on a donn\'e les bonnes d\'efinitions, on doit
obtenir que l'application $\cal F$ transforme les
param\`etres stables au sens de \ref{traduction} en des
combinaisons lin\'eaires stables \`a la Langlands; on renvoie
aux paragraphes suivants pour expliquer cette derni\`ere
notion. Cela a d\'ej\`a \'et\'e fait dans \cite{mw} par. 6
dans ce qui est, en fait, le cas le plus difficile; en effet
la difficult\'e vient de ce qu'il faut travailler avec des
param\`etres elliptiques et non pas des param\`etres discrets
et cette  difficult\'e n'appara\^{\i}t vraiment que quand
$[VP(\chi)]$ contient $+1$ et/ou $-1$. 

On dit qu'une orbite unipotente  d'un groupe lin\'eaire
complexe est elliptique symplectique (resp. orthogonale) si
ses blocs de Jordan sont tous pairs (resp. impairs)
intervenant avec multiplicit\'e au plus 2; pour le calcul du
commutant, il n'y a pas de changement majeur au lieu d'avoir
des groupes $O(1)$, on a un groupe $O(2)$ chaque fois qu'il y
a multiplicit\'e 2 (cf. \ref{classification}). Disons qu'un
param\`etre
$\psi,\epsilon$ est elliptique de niveau 0
 si $\psi_{\vert W_F}$ est mod\'er\'ement ramifi\'e comme en
\ref{classification} et avec les notations de loc.cite, pour
tout
$[u]\in [VP(\chi)]$ et pour tout $\zeta=\pm 1$, l'orbite
$U_{[u],\zeta}$ est elliptique symplectique ou orthogonale
(la diff\'erence entre symplectique et orthogonale \'etant
comme en \ref{classification}). On remarque qu'il y a une
diff\'erence pour $[u]\neq [\pm 1]$ et pour $[u]=[\pm 1]$. En
effet dans le premier cas, il revient au m\^eme de dire que
$U_{[u]}$ est discr\`ete (resp. elliptique) que de dire que
chaque $U_{[u],\zeta}$, pour
$\zeta=\pm 1$ est discr\`ete ou elliptique et de plus, ce qui
est le plus int\'eressant est que $U_{[u]}$ d\'etermine
chaque $U_{[u],\zeta}$. Ce n'est plus le cas si
$[u]=[\pm 1]$; dans ce deuxi\`eme cas la multiplicit\'e d'un
bloc de Jordan dans
$U_{[u]}$ a comme seule obligation d'\^etre inf\'erieure ou
\'egale \`a 4 et le point le plus grave est que $U_{[u]}$ ne
d\'etermine pas chaque
$U_{[u],\pm}$. 

Ici $u=\pm 1$.
On consid\`ere   l'espace vectoriel complexe de base
les quadruplets $(U_{[u],\pm},\epsilon_{[u],\pm})$ et on note
${\mathbb C}[Ell_{[u]}]$ son sous-espace vectoriel engendr\'e
par les
\'el\'ements 
$$\sum_{\epsilon_{[u],+},
\epsilon_{[u],-}}\quad\biggl(
\prod_{\begin{array}{l}\alpha_+\in Jord(U_{[u],+});
mult_+(\alpha_+)=2,\\
\alpha_-\in
Jord(U_{[u],-});
mult_-(\alpha_-)=2
\end{array}}\epsilon_{[u],+}(\alpha_+)\epsilon_{[u],-}
(\alpha_-)\biggr) (U_{[u],\pm},\epsilon_{[u],\pm})$$o\`u
$U_{[u],\pm}$ est fix\'e et la somme ne porte que sur les
$\epsilon_{[u],\pm}$ fix\'es sur l'ensemble des
$\alpha_{\pm}$ dans $Jord(U_{[u],\pm})$ dont la multiplicit\'e
$mult_{[u],\pm}(\alpha)$ est 1.

On a d\'efini en \cite{mw} 6.11 une involution de ${\mathbb
C}[Ell_{[u]}]$; il faut transporter le $\cal F$ du (i) de
loc.cit par la bijection
$rea$ du (ii) de loc. cit.. C'\'etait m\^eme une isom\'etrie,
mais on n'insiste pas la-dessus ici. C'est trop technique pour
qu'on redonne  la d\'efinition. On note
${\cal F}_{[u]}$ cette involution.

Consid\'erons maintenant le cas de $[u]\neq \pm 1$;  on pose
ici
${\mathbb C}[Disc_{[u]}]$, l'espace vectoriel complexe de
base les \'el\'ements $U_{[u]},\epsilon_{[u]}$ o\`u est
$U_{[u]}$ est discr\`ete (c'est-\`a-dire que tous ses blocs
de Jordan ont multiplicit\'e 1).  Pour d\'efinir l'application
${\cal F}_{[u]}$, on d\'efinit le produit scalaire:
\[<(U_{[u]},\epsilon_{[u]}),(U'_{[u]},\epsilon'_{[u]})
>_{[u]}:=\begin{array}{l}=0, \text{ si } U_{[u]}\neq
U'_{[u]},\\ =\sigma(U_{[u]})\sigma(\epsilon_{[u]})
\sigma_{[u]}(\epsilon'_{[u]})\prod_{\alpha\in Jord
(U_{[u]});\epsilon_{[u]}(\alpha)=-1}\epsilon'_{[u]}(\alpha)
\text{ sinon,}
\end{array}\] o\`u tous les $\sigma$ sont des signes
d\'ependant de l'objet dans la parenth\`ese;  ici on n'a
besoin que de
$\sigma_{[u]}(\epsilon'_{[u]})$. On le prend \'egal \`a
$\times_{\alpha\in Jord(U_{[u]}),\alpha\equiv
1[2]}\epsilon'_{[u]}(\alpha)$.

On pose alors
${\cal F}_{[u]}(U_{[u]},\epsilon_{[u]}):=\sum_{ 
\epsilon'_{[u]}}<(U_{[u]},\epsilon_{[u]}),
(U_{[u]},\epsilon'_{[u]})>_{[u]}(U_{[u]},\epsilon'_{[u]})$.
On remarque ais\'ement que ${\cal F}_{[u]}^2= 2^{\vert
Jord(U_{[u]})\vert} {\cal F}_{[u]}$.

Pour homog\'en\'eiser, on d\'efinit aussi  ${\mathbb
C}[Ell_{[u]}]$; c'est l'espace vectoriel engendr\'e par les
\'el\'ements:
$$\sum_{\epsilon_{[u]}} \biggl(
\prod_{\alpha\in Jord(U_{[u]});
mult(\alpha)=2}\epsilon_{[u]}(\alpha)\biggr)
(U_{[u]},\epsilon_{[u]})$$o\`u
$U_{[u],\pm}$ est fix\'e et la somme ne porte que sur les
$\epsilon_{[u],\pm}$ fix\'es sur l'ensemble des
$\alpha_{\pm}$ dans $Jord(U_{[u],\pm})$ dont la multiplicit\'e
$mult_{[u],\pm}(\alpha)$ est 1. On \'etend ${\cal F}_{[u]}$
\`a ${\mathbb C}[Ell_{[u]}]$ en \'etendant la formule
d\'ej\`a donn\'ee en pr\'ecisant simplement que le produit ne
porte que sur les
$\alpha$ dont la multiplicit\'e comme bloc de Jordan est 1.

On pose ${\mathbb C}[Ell_{\chi}]:=\otimes_{[u]\in [VP(\chi)];
[u]} {\mathbb C}[Ell_{[u]}]$. Et on d\'efinit ${\cal
F}:=\otimes_{[u]\in [VP(\chi)]} {\cal F}_{[u]}$

\subsection{Restriction aux parahoriques des
repr\'esentations\label{conjecture}} 
On a d\'efini la repr\'esentation de
Springer-Lusztig en
\ref{springerlusztig}; suivant \cite{m} 6.5, on modifie
l\'eg\`erement cette d\'efinition dans le cas des groupes
unitaires, on la note alors $SpL_{ell}$; cela induit alors un
changement dans \ref{springerlusztig} que l'on marque par le
changement de notation de $SpL_{ell}$. Soit
$(\psi,\epsilon)$ un param\`etre discret (ou elliptique) de
niveau 0;  on note
$\epsilon_{Z}$ la restriction de $\epsilon$ \`a l'\'el\'ement
non trivial du centre de
$Sp(2n,{\mathbb C})$. Pour $\sharp=iso$ ou $an$, on dit que
$\epsilon_Z=\sharp$ si $\epsilon_Z=1$ quand $\sharp=iso$ et
$-1$ sinon. On note $\vert D\vert$ l'involution de \cite{au}
et \cite{ss} qui envoie une repr\'esentation irr\'eductible
sur une repr\'esentation irr\'eductible.

\bf Conjecture: \sl Il existe  une bijection entre les
param\`etres discrets de niveau 0 ayant $\chi$ comme
restriction \`a  $I_F$ et v\'erifiant $\epsilon_Z=\sharp$ et
les s\'eries discr\`etes de niveau z\'ero du groupe
$SO(2n+1,F)_{\sharp}$ ayant
$\chi$ comme
\'el\'ement semi-simple de leur support cuspidal:
$(\psi,\epsilon)\mapsto \pi_{\psi,\epsilon}$  qui s'\'etend
en une bijection, not\'ee $rea$, entre ${\mathbb
C}[Ell_{\chi}]$ et l'espace vectoriel complexe engendr\'e par
les repr\'esentations elliptiques au sens d'Arthur  ayant
$\chi$ comme \'el\'ement semi-simple de leur support cuspidal
avec la propri\'et\'e:
 pour tout param\`etre discret de niveau z\'ero,
$(\psi,\epsilon)$, 
$k_\chi (\rho\circ
\iota )\, SpL_{ell}(\psi,\epsilon)$ est un pseudo-coefficient
de
$
\,\vert D\vert rea\, {\cal F} (\psi,\epsilon)$ (ou plus
exactement
a les m\^emes int\'egrales orbitales qu'un
pseudo-coefficient en les points semi-simples r\'eguliers
elliptiques de $SO(2n+1,F)_{\sharp}$).
\rm

 Cette conjecture est d\'emontr\'ee dans
\cite{w2} pour les repr\'esentations de r\'eduction
unipotente.

\

\bf Remarque: \sl la conjecture r\'esulte de (\cite{m} 7.)
modulo le r\'esultat annonc\'e dans \cite{auber} comme
expliqu\'e dans l'introduction; mais quand \cite{auber} sera
disponible il faudra s'assurer que le signe qui s'y
introduit compense bien la diff\'erence entre la
d\'efinition de ${\cal F}$ ici et celle de \cite{m}; en
\cite{m}, le signe qui s'introduit est
$\prod_{u\neq \pm 1}\prod_{\alpha\in
Jord(U_{[u]});\alpha\equiv m([u])+1[2]}\epsilon'(\alpha)$
alors qu'ici on a fait le produit sur les blocs de Jordan
impair.\rm

\

\subsection{interpr\'etation des r\'esultats de stabilit\'e}

On a vu en \cite{mw} 4.6 (\`a la suite d'Arthur) que la
stabilit\'e des repr\'esentations elliptiques se lit sur les
int\'egrales orbitales des pseudo-coefficients, en admettant
la conjecture de
\ref{conjecture} on peut d\'ecrire les combinaisons
lin\'eaires de repr\'esentations discr\`etes qui sont
stables. Soit $\psi,\epsilon$ un param\`etre discret de
niveau 0 de restriction le caract\`ere $\chi$ \`a $I_F$.  On
note encore
$\epsilon_Z$ la restriction de $\epsilon$ au centre de
$SO(2n+1,F)_{\sharp}$ (o\`u $\sharp=iso$ ou $an$) et on dit
que
$\epsilon_Z=\sharp$ si
$\epsilon_Z=+$ quand $\sharp=iso$ et $-$ quand $\sharp=an$.

\bf Th\'eor\`eme: \sl Ici on admet la conjecture de
\ref{conjecture}. Soient $\sharp=iso$ ou $an$ et  $\psi$ un
param\`etre discret de niveau 0.

(i) La combinaison lin\'eaire:
$$
\sum_{\epsilon; \epsilon_Z=\sharp}
\pi_{\psi,\epsilon}$$  est stable pour le groupe
$SO(2n+1,F)_{\sharp}$. De plus dans le transfert entre
$SO(2n+1,F)_{an}$ et $SO(2n+1,F)_{iso}$, les combinaisons
lin\'eaires
$\epsilon_Z\sum_{\epsilon;
\epsilon_Z=\sharp}
\,
\pi_{\psi,\epsilon}$ se correspondent (ici $\sharp$ est vu
comme un \'el\'ement de $\pm 1$).

(ii) toute combinaison lin\'eaire des repr\'esentations
$\pi_{\psi,\epsilon}$ pour $\epsilon$ variant avec
$\epsilon_{Z}=\sharp$ est instable si elle n'est pas
proportionnelle \`a la combinaison \'ecrite en (i).

\rm On fixe $\sharp=iso$ ou $an$ et
 on note
${\mathbb C}[Ell_{\chi}]_{stable,\sharp}$ le sous-espace de
${\mathbb C}[Ell_{\chi}]$ form\'e de l'image par $rea$ des
combinaisons lin\'eaires stables de repr\'esentations
elliptiques pour $SO(2n+1,F)_{\sharp}$. Et on note ${\mathbb
C}[Ell_{\chi}]_{st,sst}$ le sous-espace de ${\mathbb
C}[Ell_{\chi}]$ engendr\'e par les param\`etres elliptiques de
niveau 0, stables ou semi-stables. On reprend les notations
 $U'([u])$ et $U''([u])$ de \ref{springerlusztig},  l'espace
ci-dessus est donc naturellement la somme directe des 2
sous-espaces, l'un (resp. l'autre) engendr\'e par les
param\`etres
$(\psi,\epsilon)$ tels que
$U''_{[u]}=0$ (resp. $U'_{[u]}=0$)  pour tout
$[u]\in [VP(\chi)]$. Avec la conjecture et le th\'eor\`eme de
\ref{theoreme} (compl\'et\'e par la remarque de
\ref{stabilite}), on sait que
$rea\, \circ\, {\cal F}$ induit entre ${\mathbb
C}[Ell_{\chi}]_{st,sst}$ et
${\mathbb C}[Ell_{\chi}]_{stable,iso}\oplus {\mathbb
C}[Ell_{\chi}]_{stable,an}$. Il suffit donc de calculer
l'image par ${\cal F}$ d'un param\`etre
$(\psi,\epsilon)$ elliptique de niveau 0  qui soit stable ou
semi-stable et de reprojeter sur l'espace vectoriel
engendr\'e par les param\`etres $(\psi,\epsilon)$ v\'erifiant
$\epsilon_Z=\sharp$ quand $\sharp$ est fix\'e. Fixons donc
$\zeta=\pm$ et calculons l'image
${\cal F}(\psi,\epsilon)$ en supposant que pour tout $[u]\in
[VP(\chi)]$ l'orbite $U^\delta_{[u]}=0$, o\`u $\delta='$ si
$\zeta=+$ et $\delta=''$ si $\zeta=-$. On esp\`ere que le
lecteur comprendra une d\'ecomposition
${\cal F}(\psi,\epsilon)=\times_{[u]\in [VP(\chi)]}{\cal
F}_{[u]}(\psi_{[u]},\epsilon_{[u]})$. Et on calcule ${\cal
F}_{[u]}(\psi_{[u]},\epsilon_{[u]})$ en supposant d'abord 
que $[u]\neq [\pm 1]$; d'apr\`es la d\'efinition, on a 
\[{\cal
F}_{[u]}(\psi_{[u]},\epsilon_{[u]})=\sigma(U_{[u]})\sigma_{[u]}(\epsilon_{[u]}
)\sum_{\epsilon'_{[u]}}
\sigma_{[u]}(\epsilon'_{[u]})\biggl(
\prod_{\begin{array}{l}\alpha\in Jord(U_{[u]});\\
\epsilon_{[u]}(\alpha)=-1
\end{array}}\epsilon'_{[u]}(\alpha)\biggr)
(\psi_{[u]},\epsilon'_{[u]}).\] Or
$\epsilon_{[u]}(\alpha)=\zeta$ par hypoth\`ese pour tout
$\alpha$; la formule ci-dessus se simplifie donc si $\zeta=+$
en 
\[{\cal F}_{[u]}(\psi_{[u]},\epsilon_{[u]})=
\sigma_{[u]}(\epsilon_{[u]})
\sum_{\epsilon'_{[u]}}
\sigma_{[u]}(\epsilon'_{[u]}) (\psi_{[u]},\epsilon'_{[u]}).\] 
 Par contre si $\zeta=-$, elle se simplifie en:
\[{\cal F}_{[u]}(\psi_{[u]},\epsilon_{[u]})=
\sigma_{[u]}(\epsilon_{[u]})
\sum_{\epsilon'_{[u]}}
\sigma_{[u]}(\epsilon'_{[u]} )\biggl(
\prod_{\alpha\in Jord(U_{[u]})}\epsilon'_{[u]}(\alpha)\biggr)
(\psi_{[u]},\epsilon'_{[u]}).\]  Le cas de $[u]=[\pm 1]$ est
exactement celui trait\'e en
\cite{mw} 6.12, et le r\'esultat est analogue \`a ci-dessus.

En revenant au produit, on obtient dans le cas $\zeta=+$,
avec $\sigma$ un signe qui d\'epend de $\psi,\epsilon$:
\[{\cal
F}(\psi,\epsilon)=\sigma\sum_{\epsilon'}
 (\psi,\epsilon').\] Dans le cas
$\zeta=-$, dans la formule s'ajoute
$\prod_{\alpha\in \cup_{[u]}Jord(U_{[u]})}\epsilon'(\alpha)$
qui n'est autre que $\epsilon'_Z$. Quand on s\'epare les 2
morceaux, celui correspondant \`a $\sharp=iso$ et $\sharp=an$,
$\epsilon'_Z$ est constant dans chaque morceau.

Pour pouvoir en d\'eduire le r\'esultat de stabilit\'e
cherch\'ee, il faut utiliser la conjecture \ref{conjecture}
qui permet de calculer les int\'egrales orbitales des
caract\`eres des repr\'esentations pour les \'el\'ements
elliptiques. Mais il faut d'abord enlever $\vert D\vert$. Or
pour toute repr\'esentation $\pi$ irr\'eductible et pour tout
\'el\'ement elliptique $\gamma$ de $G$ le caract\`ere de
$\pi$ et de $\vert D\vert \pi$ co\"{\i}ncident en $\gamma$ au
signe $(-1)^{rg_G-rg_{P_{cusp,\pi}}}$, o\`u $P_{cusp,\pi}$
est le sous-groupe parabolique de $G$ minimal pour la
propri\'et\'e que $res_P(\pi)$ est non nulle. Pour $\pi$ de
la forme $\pi(\psi,\epsilon)$,
$(-1)^{rg_G-rg_{P_{cusp,\pi}}}=\prod_{\alpha\in
\cup_{[u]}Jord(U_{[u]});\alpha\equiv 0[2]}\epsilon(\alpha)$.
On a ainsi d\'emontr\'e que la distribution
\[\sum_{\epsilon'}\biggl(\prod_{\alpha\in
\cup_{[u]}Jord(U_{[u]})}\epsilon'(\alpha)\biggr)\pi(\psi,\epsilon')\]
est stable et que pour $\sharp=iso$ ou $an$ les seules
distributions stables sont les sous-sommes de la somme
ci-dessus o\`u l'on ne somme que sur les $\epsilon'$ tels que
$\epsilon'_{Z}=\sharp$. Cela donne le r\'esultat annonc\'e.

\end{document}